\documentclass[11pt]{article}

\usepackage{indentfirst}
\usepackage{amsfonts}
\usepackage{amsmath}
\usepackage{amssymb}
\usepackage{amsbsy, amsthm}

\newtheorem{dfn}{Definition} [section]
\newtheorem{theorem}[dfn]{Theorem}
\newtheorem{lemma}[dfn]{Lemma}

\newtheorem{conjecture}[dfn]{Conjecture}
\newenvironment{pf}{\noindent{\bf Proof.}}
{\enspace\vrule height5pt depth0pt width5pt}

\usepackage{latexsym,amsmath,amssymb,amsfonts,epsfig,graphicx,cite,psfrag}
\usepackage{eepic,color,colordvi,amscd}

\newif\iftrackchanges
\newif\iftracknew
\trackchangestrue   
\tracknewtrue    
\iftrackchanges
  \usepackage{color}
  \usepackage[normalem]{ulem}

  \newcommand{\added}[2][]{{\color{blue}#2\ifx&#1&\else%
    \textsuperscript{(#1)}\fi}}
  \newcommand{\deleted}[2][]{{\color{red}\sout{#2}\ifx&#1&\else%
    \textsuperscript{(#1)}\fi}}
  \newcommand{\deletedl}[1]{{\color{red}#1}}
\else
  \iftracknew
    \newcommand{\added}[2][]{#2}
    \newcommand{\deleted}[2][]{}
    \newcommand{\deletedl}[1]{}
  \else
    \newcommand{\added}[2][]{}
    \newcommand{\deleted}[2][]{#2}
    \newcommand{\deletedl}[1]{#1}
  \fi
\fi

\begin{document}

\title{Cycle lengths and minimum degree of graphs}
\author{Chun-Hung Liu\thanks{Department of Mathematics,
Princeton University,
Princeton, New Jersey 08544, USA. Email: chliu@math.princeton.edu.}
\and
Jie Ma\thanks{School of Mathematical Sciences, University of Science and Technology of China, Hefei, Anhui 230026, China. Email: jiema@ustc.edu.cn. Partially supported by NSFC project 11501539.}
}

\date{}

\maketitle

\begin{abstract}
\noindent
There has been extensive research on cycle lengths in graphs with large minimum degree.
In this paper, we obtain several new and tight results in this area.
Let $G$ be a graph with minimum degree at least $k+1$.
We prove that if $G$ is bipartite, then there are $k$ cycles in $G$ whose lengths form an arithmetic progression with common difference two.
For general graph $G$, we show that $G$ contains $\lfloor k/2\rfloor$ cycles with consecutive even lengths and $k-3$ cycles whose lengths form an arithmetic progression with common difference one or two.
In addition, if $G$ is 2-connected and non-bipartite, then $G$ contains $\lfloor k/2\rfloor$ cycles with consecutive odd lengths.

Thomassen (1983) made two conjectures on cycle lengths modulo a fixed integer $k$: (1) every graph with minimum degree at least $k+1$ contains cycles of all even lengths modulo $k$; (2) every 2-connected non-bipartite graph with minimum degree at least $k+1$ contains cycles of all lengths modulo $k$.
These two conjectures, if true, are best possible.
Our results confirm both conjectures when $k$ is even.
And when $k$ is odd, we show that minimum degree at least $k+4$ suffices.
This improves all previous results in this direction.
Moreover, our results derive new upper bounds of the chromatic number in terms of the longest sequence of cycles with consecutive (even or odd) lengths.
\end{abstract}

\section{Introduction}
\noindent 
%
The study of the distribution of cycle lengths is a fundamental area in modern graph theory, which has led to numerous results in abundant subjects.
A common practice is investigating if certain graph properties, such as large average degree, large chromatic number, large connectivity, or nice expansion properties, are sufficient to ensure the existence of cycles of some particular lengths.
In this article, all graphs are simple and we consider the distribution of cycle lengths in graphs with large minimum degree, aiming to understand the relation between cycle lengths and minimum degree in great depth.

One classical result in this direction is due to Dirac \cite{Dir52} in 1950s: every graph $G$ with $n\geq 3$ vertices and with minimum degree at least $n/2$ contains a {\em Hamilton} cycle (i.e., a cycle passing through all vertices of $G$).
Since then, there has been extensive research to investigate cycle lengths in graphs $G$ with large minimum degree $\delta(G)$, where $\delta(G)$ desponds on $|V(G)|$.
To name a few, \cite{Alon86,BFG98,BH89} are about the length of the longest cycle, \cite{Hagg} is about the existence of cycles with specified lengths, and \cite{Bon71,EFGS,BS74,GHS02,NS04,NS06} are about the range of cycle lengths.

However, it is more general if the minimum degree is independent with the number of vertices.
Dirac \cite{Dir52} proved that every 2-connected graph with $n$ vertices and minimum degree $k$ contains a cycle of length at least $\min\{n,2k\}$.
Voss and Zuluaga \cite{VZ77} generalized this by proving that every 2-connected non-bipartite graph with $n$ vertices and  minimum degree $k$ contains an even cycle of
length at least $\min\{n, 2k\}$ and an odd cycle of length at least $\min\{n, 2k-1\}$.
Bondy and Vince \cite{BV98} solved a question of Erd\H{o}s by proving 
that if all but at most two vertices of $G$ have degree at least three, then there are two cycles in $G$ whose lengths differ by one or two.
H\"{a}ggkvist and Scott \cite{HS} proved that every connected cubic graph other than $K_4$ contains two cycles whose lengths differ by two.

Bondy and Vince's theorem was improved by several authors.
H\"{a}ggkvist and Scott \cite{HS98} proved that every graph with minimum degree $\Omega(k^2)$ contains $k$ cycles of consecutive even lengths.
Verstra\"ete \cite{V00} improved this quadratic bound to be linear by proving that every graph with average degree at least $8k$ and even girth $g$ contains $(g/2-1)k$ cycles of consecutive even lengths.
In \cite{SV08}, Sudakov and Verstra\"ete further pushed the number of lengths of the cycles to be exponential: every graph with average degree $192(k+1)$ and girth $g$ contains $k^{\lfloor(g-1)/2\rfloor}$ cycles of consecutive even lengths.
Very recently, the second author \cite{Ma} obtained an analogue for odd cycle: every 2-connected non-bipartite graph with average degree $456k$ and girth $g$ contains $k^{\lfloor(g-1)/2\rfloor}$ cycles of consecutive odd lengths.
On the other hand, without considering the parity of the cycles, Fan \cite{Fan02} obtained similar results with better minimum degree conditions by proving the following result.
Every graph $G$ with minimum degree $\delta(G)\ge 3k$ contains $k+1$ cycles $C_0,C_1,...,C_k$ such that $|E(C_0)|>k+1, |E(C_{i})|-|E(C_{i-1})|=2$ for all $1\le i\le k-1$ and $1\le |E(C_{k})|-|E(C_{k-1})|\le 2$, and furthermore, if $\delta(G)\ge 3k+1$, then $|E(C_{k})|-|E(C_{k-1})|=2$.
In the same paper \cite{Fan02}, he also resolved a problem of Bondy and Vince \cite{BV98} by showing that every 3-connected non-bipartite graph $G$ with $\delta(G)\ge 3k$ contains $2k$ cycles with consecutive lengths $m,m+1,...,m+2k-1$ for some integer $m\ge k+2$.

To better understand the above results, we remark that in order to ensure two or more odd cycle lengths, 2-connectedness is necessary in addition to the non-bipartiteness.
There exist infinitely many non-bipartite connected graphs with arbitrary large minimum degree but containing a unique odd cycle:
for arbitrary $t$ and odd $s$, let $G$ be obtained from $s$ disjoint copies of $K_{t,t}$ and an odd cycle $C_s$ such that each $K_{t,t}$ intersects $C_s$ in exactly one vertex.

\subsection{Paths and cycles of consecutive lengths}\label{subsection 1.1}
\noindent Throughout the rest of this paper, $k$ is a fixed positive integer, unless otherwise specified.
We say that a sequence of paths or cycles $H_1,H_2,...,H_k$ satisfies the {\em length condition}
if $\lvert E(H_1) \rvert \ge 2$ and $\lvert E(H_{i+1}) \rvert - \lvert E(H_i) \rvert = 2$ for $1 \leq i \leq k-1$.
We also say that $k$ paths or $k$ cycles satisfy the length condition if they can form such a sequence.

In order to study cycles of consecutive (even or odd) lengths in graphs, we begin by considering paths in bipartite graphs.
Our first theorem says that there exist optimal number of paths in bipartite graphs between two fixed vertices and satisfying the length condition.

\begin{theorem} \label{thm:bipartite_main}
Let $G$ be a 2-connected bipartite graph and $x,y$ distinct vertices of $G$.
If every vertex in $G$ other than $x,y$ has degree at least $k+1$, then there exist $k$ paths $P_1,P_2,...,P_k$ from $x$ to $y$ in $G$ with the length condition.
\end{theorem}

\noindent We point out that this result is crucial to the proofs of all other results in this paper.
The minimum degree condition in Theorem \ref{thm:bipartite_main} is tight for infinitely many graphs, by considering the complete bipartite graphs $K_{k,n}$ for all $n\ge k$, where $x,y$ are two vertices in the part of size $k$.

\medskip

The following theorem on cycles in bipartite graphs can be derived from Theorem \ref{thm:bipartite_main}.

\begin{theorem} \label{intro:cycle bipartite}
Let $G$ be a bipartite graph and $v$ a vertex of $G$.
If every vertex of $G$ other than $v$ has degree at least $k+1$, then $G$ contains $k$ cycles with the length condition.
\end{theorem}

\noindent An immediate corollary of Theorem \ref{intro:cycle bipartite} is that every bipartite graph with minimum degree at least $k+1$ contains $k$ cycles with the length condition.
The complete bipartite graphs $K_{k,n}$ for all $n\ge k$ also show the tightness of the minimum degree condition.

\medskip

We then investigate cycle lengths in general graphs.

\begin{theorem} \label{intro:even/odd cycle in general graph}
If the minimum degree of graph $G$ is at least $k+1$, then $G$ contains $\lfloor k/2\rfloor$ cycles with consecutive even lengths.
Furthermore, if $G$ is 2-connected and non-bipartite, then $G$ contains $\lfloor k/2\rfloor$ cycles with consecutive odd lengths.
\end{theorem}

\noindent
We see that Theorem \ref{intro:even/odd cycle in general graph} is tight, as the complete graph $K_{k+2}$ has exactly $\lfloor k/2\rfloor$ different even cycle lengths regardless of the parity of $k$, and it has exactly $\lfloor k/2\rfloor$ different odd cycle lengths when $k$ is even.

In the coming two theorems, we consider 3-connected and 2-connected non-bipartite graphs respectively.

\begin{theorem} \label{intro:3-connected non-bipartite}
If $G$ is a 3-connected non-bipartite graph with minimum degree at least $k+1$, then $G$ contains $2\lfloor \frac{k-1}{2}\rfloor$ cycles with consecutive lengths.
\end{theorem}

\begin{theorem}\label{intro:2-con non-bipa cycles}
If $G$ is a 2-connected non-bipartite graph with minimum degree at least $k+3$, then $G$ contains $k$ cycles with consecutive lengths or the length condition.
\end{theorem}

\noindent Theorem \ref{intro:3-connected non-bipartite} improves a result of Fan \cite{Fan02}, which was originally asked by Bondy and Vince \cite{BV98}.
Note that Bondy and Vince \cite{BV98} constructed an infinite family of 2-connected non-bipartite graphs with arbitrarily large minimum degree but containing no two cycles whose lengths differ by one.
So the connectivity condition in Theorem \ref{intro:3-connected non-bipartite} cannot be lowered, and the conclusion for cycles with the length condition in Theorem \ref{intro:2-con non-bipa cycles} cannot be dropped.
Moreover, every graph on at most $2k$ vertices does not have $k$ cycle with the length condition.
Hence, $K_{2k}$ is an example showing that the conclusion for cycles with consecutive lengths in Theorem \ref{intro:2-con non-bipa cycles} also cannot be removed when $k \geq 4$.
(But Theorem \ref{intro:even/odd cycle in general graph} ensures the existence of cycles with the length condition when $k=2$.)
Therefore, Theorem \ref{intro:2-con non-bipa cycles} cannot be further improved to require only cycles with consecutive lengths or only cycles with the length condition in general.
By considering complete graphs of certain orders, we can see that the difference between the minimum degree conditions in Theorems \ref{intro:3-connected non-bipartite} and \ref{intro:2-con non-bipa cycles} and the optimal bounds is at most two.

The next result studies cycle lengths in general graphs, without assuming connectivity and bipartiteness.

\begin{theorem}\label{intro:more general graph k cycles}
If $G$ is a graph with minimum degree at least $k+4$,
then $G$ contains $k$ cycles with consecutive lengths or the length condition.
\end{theorem}

\noindent This improves some aforementioned results in \cite{V00,Fan02}.
We direct readers to Section \ref{sec:remarks} for a discussion on the tightness of this theorem.

\subsection{Cycle lengths modulo $k$}\label{subsection 1.2}
\noindent The study of cycle lengths modulo an integer $k$ can be dated to Burr and Erd\H{o}s (See \cite{Erd76}).
They conjectured that there exists a constant $c_k$ for each odd $k$ such that every graph with average degree at least $c_k$ contains cycles of all lengths modulo $k$.
This conjecture was resolved by Bollob\'as in \cite{B77}, where he proved that $c_k\le 2[(k+1)^k-1]/k$.
Thomassen \cite{Th83,Th88} generalized this by showing that
every graph $G$ with minimum degree at least $4k(k+1)$ contains cycles of all lengths $m$ modulo $k$, except when $m$ is odd and $k$ is even.
Note that the exceptional case is needed, as when $k$ is even and $G$ is bipartite, there is no odd cycle in $G$ and thus no cycle of odd length $m$ modulo $k$.
Thomassen \cite{Th83} observed that $K_{k+1}$ has no cycle of length 2 modulo $k$, and made the following conjecture.

\begin{conjecture}[Thomassen \cite{Th83}]\label{conj:Thom1}
For every positive integer $k$, every graph with minimum degree at least $k+1$ contains cycles of all even lengths modulo $k$.
\end{conjecture}

Thomassen \cite{Th83} also proved that there exists a function $\theta(k)$ for every $k$ such that every 2-connected non-bipartite graph with minimum degree at least $\theta(k)$ contains cycles of all lengths modulo $k$.
Note that the same graphs defined before Section \ref{subsection 1.1} show that 2-connectivity and non-bipartiteness are necessary conditions here (for even $k$).

\begin{conjecture}[Thomassen \cite{Th83}]\label{conj:Thom2}
For every positive integer $k$, every 2-connected non-bipartite graph with minimum degree at least $k+1$ contains cycles of all lengths modulo $k$.\footnote{It is quoted from \cite{Th83} that ``$K_{k+2}$ shows that $\theta(k)\ge k+2$. It is tempting to conjecture that equality holds." Since $K_{k+2}$ does contain cycles of all lengths modulo $k$, we believe that it meant to conjecture $\theta(k)=k+1$.}
\end{conjecture}

It is known that the minimum degree $\Omega(k)$ suffices for both Conjectures \ref{conj:Thom1} and \ref{conj:Thom2}.
A theorem of Verstra\"ete \cite{V00} implies that for all $k$, every graph with average degree at least $8k$ contains cycles of all even lengths modulo $k$.
For all odd $k$, a result of Fan \cite{Fan02} shows that minimum degree at least $3k-2$ suffices.
Diwan \cite{D10} obtained a better bound for Conjecture \ref{conj:Thom1}
that for every positive integer $k$, every graph $G$ with minimum degree at least $2k-1$ contains cycles of all even lengths modulo $k$, and every graph with minimum degree at least $k+1$ contains a cycle of length 4 modulo $k$.
For Conjecture \ref{conj:Thom2}, a recent result of \cite{Ma} about consecutive odd cycles implies that minimum degree $\Omega(k)$ is suffices to ensure the existence of cycles of all lengths modulo $k$.

Using our results in Section \ref{subsection 1.1}, we obtain several consequences on cycle lengths modulo $k$, which improve all previous bounds on Conjectures \ref{conj:Thom1} and \ref{conj:Thom2}.
In particular, the following theorem settles both Conjectures \ref{conj:Thom1} and \ref{conj:Thom2} for all even integers $k$.

\begin{theorem}\label{intro:cycles mod even k}
Let $k$ be a positive even integer.
If $G$ is a graph with minimum degree at least $k+1$, then $G$ contains cycles of all even lengths modulo $k$.
Furthermore, if $G$ is 2-connected and non-bipartite, then $G$ contains cycles of all lengths modulo $k$.
\end{theorem}

The case for odd $k$ seems more intricate than the case for even $k$.
The next two theorems can be derived from Theorems \ref{intro:2-con non-bipa cycles} and \ref{intro:more general graph k cycles}, respectively.

\begin{theorem}\label{intro:2-con non-bipa cycles modulo k}
Let $k$ be a positive odd integer. If $G$ is a 2-connected non-bipartite graph with minimum degree at least $k+3$,
then $G$ contains cycles of all lengths modulo $k$.
\end{theorem}

\begin{theorem}\label{intro:general cycles mod k}
Let $k$ be a positive odd integer. If $G$ is a graph with minimum degree at least $k+4$,
then $G$ contains cycles of all lengths modulo $k$.
\end{theorem}

\noindent In other words, when $k$ is odd, the difference between the minimum degree conditions of our results and the bounds of Thomassen's conjectures is at most three.

\subsection{Cycles of consecutive lengths and chromatic number}\label{subsection 1.3}
\noindent
The chromatic number and the length of cycles are also related.
Diwan, Kenkre and Vishwanathan \cite{DKV} conjectured that for every pair of integers $m$ and $k$, if graph $G$ has no cycle of length $m$ modulo $k$, then the chromatic number of $G$ is at most $k+o(k)$.
This was resolved by Chen, Ma and Zang in a recent paper\cite{CMZ}, where they also studied the relations between cycle lengths modulo $k$ and chromatic number of digraphs.

Given a graph $G$, define $L_e(G)$ and $L_o(G)$ to be the sets of even and odd cycle lengths in $G$, respectively.
We define $ce(G)$ and $co(G)$ to be the largest integers $m$ and $n$, respectively,
such that $G$ contains $m$ cycles of consecutive even lengths and $n$ cycles of consecutive odd lengths.
And we denote the largest integer $\ell$ by $c(G)$ such that $G$ contains $\ell$ cycles of consecutive lengths.

We say that a graph $G$ is {\em $k$-chromatic} if its chromatic number $\chi(G)$ equals $k$.
It is well-known that every $k$-chromatic graph has a cycle of length at least $k$. In 1966, Erd\H{o}s and Hajnal \cite{EH66} provided an analogue that every $k$-chromatic graph has an odd cycle of length at least $k-1$.
Confirming a conjecture of Bollob\'as and Erd\H{o}s, Gyarf\'as \cite{Gyarfas} generalized the result of Erd\H{o}s and Hajnal by showing that every graph $G$ satisfies $\chi(G)\le 2|L_o(G)|+2$. Mihok and Schiermeyer \cite{MS04} proved that $\chi(G)\le 2|L_e(G)|+3$ for every graph $G$.
Recently, Kostochka, Sudakov and Verstra\"ete \cite{KSV} proved a conjecture of Erd\H{o}s \cite{Erd92} that every triangle-free $k$-chromatic graph $G$ contains at least $\Omega(k^2\log k)$ cycles of consecutive lengths.

Using Theorem \ref{intro:even/odd cycle in general graph}, we obtain a new upper bound of the chromatic number in terms of the longest sequence of consecutive even or odd cycle lengths.

\begin{theorem}\label{intro:ce(G)+co(G)}
For every graph $G$, $\chi(G)\le 2 \min\{ce(G),co(G)\}+3$.
\end{theorem}

\noindent This strengthens the result of Mihok and Schiermeyer \cite{MS04}, as clearly $ce(G)\le |L_e(G)|$.
In addition, Theorem \ref{intro:ce(G)+co(G)} is tight for the complete graphs on odd number of vertices, as $\min\{ce(K_{2k+3}),co(K_{2k+3})\}=k$.

Moreover, we show that the chromatic number can be bounded from above by the longest sequence of consecutive cycle lengths.

\begin{theorem}\label{intro:c(G)}
For every graphs $G$, $\chi(G)\le c(G)+4$.
\end{theorem}

\noindent On the other hand, complete graphs show that $\chi(G)\ge c(G)+2$.

\subsection{Notation and organization}
\noindent Let $G$ be a graph and $X$ a subset of $V(G)$.
We denote the set of vertices not in $X$ but adjacent to some vertex in $X$ by $N_G(X)$, and we define $N_G[X]:= N_G(X) \cup X$.
If $X=\{x\}$, we simply write $N_G(x)$ and $N_G[x]$ instead.
For a subgraph $D$ of $G$, we define $N_G(D):= N_G(V(D))$ and $N_G[D]:=N_G[V(D)]$.
Often we drop the subscript when $G$ is clear from context.
For a vertex $v$ of $G$, the {\em degree} of $v$, denoted by $d_G(v)$, is the number of edges in $G$ incident with $v$, and we define $d_X(v):= \lvert N_G(v)\cap X \rvert$.
A vertex is a {\em leaf} in $G$ if it has degree one in $G$.
For $S \subseteq V(G)$, we denote the subgraph of $G$ induced on $V(G)-S$ by $G-S$; for $S \subseteq E(G)$, we denote the graph $(V(G),E(G)-S)$ by $G-S$.
When $S \subseteq V(G) \cup E(G)$ with $\lvert S \rvert=1$, we write $G-S$ as $G-s$, where $s$ is the unique element of $S$.
When we identify a subset $S$ of $V(G)$, we always delete all resulting loops and parallel edges to keep the graph simple.

A pair $(A,B)$ of subsets of $V(G)$ is a {\em separation} of $G$ of {\em order $k$},
if $V(G)=A\cup B$, $|A\cap B|=k$ and $G$ has no edge with one end in $A-B$ and the other in $B-A$.
A vertex $v$ of a graph $G$ is a {\em cut-vertex} if $G-v$ contains more components than $G$.
A {\em block} $B$ in $G$ is a maximal connected subgraph of $G$ such that there exists no cut-vertex of $B$.
So a block is an isolated vertex, an edge or a 2-connected graph.
An {\em end-block} in $G$ is a block in $G$ containing at most one cut-vertex of $G$.
If $D$ is an end-block of $G$ and a vertex $x$ is the only cut-vertex of $G$ with $x \in V(B)$, then we say that $D$ is an {\em end-block with cut-vertex $x$}.
Let ${\mathcal B}(G)$ be the set of blocks in $G$ and ${\mathcal C}(G)$ be the set of cut-vertices of $G$.
The {\em block structure} of $G$ is the bipartite graph with bipartition $({\mathcal B}(G),{\mathcal C}(G))$, where $x\in {\mathcal C}(G)$ is adjacent to $B\in {\mathcal B}(G)$ if and only if $x\in V(B)$.
Note that the block structure of any graph $G$ is a forest, and it is connected if and only if $G$ is connected.
For every positive integer $k$, we say that a graph $G$ is {\em $k$-critical} if it has chromatic number $k$ and every proper subgraph of $G$ has chromatic number less than $k$.

\medskip

The rest of this paper is organized as follows.
In Section \ref{sec:path_bipartite}, we consider paths in bipartite graphs and prove Theorem \ref{thm:bipartite_main} by induction.
We then apply Theorem \ref{thm:bipartite_main} in Section \ref{sec:path_general} to 
obtain results about paths in general graphs, which will be heavily used later.
In Section \ref{sec:cycle}, we focus on cycles with the length condition and prove Theorems \ref{intro:cycle bipartite} and \ref{intro:even/odd cycle in general graph},
from which we also derive Theorems \ref{intro:cycles mod even k} and \ref{intro:ce(G)+co(G)}.
In Section \ref{sec:consecutive_cycles},
we first prove Theorem \ref{nonseparating cycle->consecutive cycles} on cycles of consecutive lengths, and then show how to derive the rest theorems mentioned in this section. 
Finally, we close the paper by mentioning some concluding remarks and open problems in Section \ref{sec:remarks}.

\section{Consecutive paths in bipartite graphs} \label{sec:path_bipartite}
\noindent We shall prove Theorem \ref{thm:bipartite_main} in this section.
To simplify the arguments, we shall prove a more general (but indeed equivalent) result. 
For this purpose, we introduce the following important concepts.
We say that $(G,x,y)$ is a {\em rooted graph} if $G$ is a graph and $x,y$ are distinct vertices of $G$. The vertices $x,y$ are called the {\em roots} of $(G,x,y)$.
A rooted graph $(G,x,y)$ is {\em bipartite} if and only if $G$ is bipartite.
The {\em minimum degree} of $(G,x,y)$ is $\min\{d_G(u): u\in V(G)-\{x,y\}\}$.
We say that $(G,x,y)$ is {\em 2-connected} if
	\begin{itemize}
		\item $G$ is a connected graph with $|V(G)|\geq 3$, and
        	\item every end-block of $G$ contains at least one of $x,y$ as a non-cut-vertex.
	\end{itemize}
Note that the block structure of $G$ is a path if $(G,x,y)$ is 2-connected.
And $x,y$ are in the same block of $G$ if and only if $G$ is 2-connected.

On the other hand, if $G$ is 2-connected, then $(G,x,y)$ is 2-connected for every pair of distinct vertices $x,y$.
Therefore, Theorem \ref{thm:bipartite_main} is an immediate corollary of the following theorem.

\begin{theorem} \label{main bipartite}
Let $(G,x,y)$ be a 2-connected bipartite rooted graph.
For any positive integer $k$, if the minimum degree of $(G,x,y)$ is at least $k+1$, then there exist $k$ paths in $G$ from $x$ to $y$ satisfying the length condition.
\end{theorem}

We shall prove Theorem \ref{main bipartite} by induction on $\lvert V(G) \rvert + \lvert E(G) \rvert$.
In the rest of this section, we define $(G,x,y)$ to be a minimum counterexample (with respect to $\lvert V(G) \rvert + \lvert E(G) \rvert$).
That is, $(G,x,y)$ is a 2-connected bipartite rooted graph with minimum degree at least $k+1$ such that $G$ does not contain $k$ paths from $x$ to $y$ satisfying the length condition; however, for any 2-connected bipartite rooted graph $(H,u,v)$ with $\lvert V(H) \rvert + \lvert E(H) \rvert< \lvert V(G) \rvert + \lvert E(G) \rvert$ and for any positive integer $r$, if the minimum degree of $(H,u,v)$ is at least $r+1$, then there are $r$ paths in $H$ from $u$ to $v$ satisfying the length condition.
By symmetry, we assume that
\begin{align}
\label{equ:dx<dy}
d_G(x)\le d_G(y).
\end{align}

Throughout the rest of this section, we will exploit related properties of $G$ and prove a series of lemmas, which will lead to the final contradiction and thus complete the proof of Theorem \ref{main bipartite}. We start by proving the following useful lemma.

\begin{lemma} \label{basis}
$|V(G)|\geq 4$, $G$ is 2-connected, and $k \geq 3$.
\end{lemma}

\begin{pf}
If $\lvert V(G) \rvert =3$, then $(G,x,y)$ has minimum degree two, so $k=1$ and the theorem follows.
Hence $\lvert V(G) \rvert \ge 4$.

Suppose that $G$ is not 2-connected.
Then there exist a cut-vertex $b$ and two connected subgraphs $G_1,G_2$ of $G$
such that $G=G_1\cup G_2$ and $V(G_1)\cap V(G_2)=\{b\}$, where $x\in V(G_1)-b$ and $y\in V(G_2)-b$.
Since $\lvert V(G_1) \rvert+\lvert V(G_2) \rvert=\lvert V(G) \rvert+1\geq 5$,
by symmetry we may assume that $\lvert V(G_1) \rvert\geq 3$.
So $(G_1,x,b)$ is 2-connected bipartite with minimum degree at least $k+1$.
By induction, there exist $k$ paths $P_1,...,P_k$ in $G_1$ from $x$ to $b$ with the length condition.
Let $P$ be a path in $G_2$ from $b$ to $y$.
Concatenating $P$ with each $P_i$ leads to $k$ paths in $G$ from $x$ to $y$ with the length condition, a contradiction. Therefore $G$ is 2-connected.

Since $G$ is 2-connected, Theorem \ref{main bipartite} is obvious when $k=1$.
The case $k=2$ can be derived by the following special case of
\cite[Corollary 3.1]{Fan02}: if $H$ is a 2-connected (not necessarily bipartite) graph and every vertex of $H$ other than two distinct vertices $u,v$ has degree at least three, then $H$ contains two paths $R_1,R_2$ from $u$ to $v$ such that $\lvert E(R_1) \rvert \geq 2$ and $1\leq \lvert E(R_2) \rvert-\lvert E(R_1)\rvert\leq 2$.
To see the implication for the case $k=2$, just notice that $G$ is bipartite and thus all paths in $G$ from $x$ to $y$ are of the same parity, implying $\lvert E(R_2) \rvert-\lvert E(R_1)\rvert= 2$.
This shows that $k\geq 3$.
\end{pf}

\begin{lemma} \label{lem:xy_notadj}
$x$ and $y$ are not adjacent in $G$.
\end{lemma}

\begin{pf}
Suppose that $x$ is adjacent to $y$ in $G$.
Let $G'=G-xy$.
Since $G$ is 2-connected, every end-block of $G'$ contains at least one of $x,y$ as non-cut-vertex.
Therefore, $(G',x,y)$ is 2-connected bipartite with minimum degree at least $k+1$.
The induction hypothesis implies that $G'$, and hence $G$, contains $k$ paths from $x$ to $y$ with the length condition, a contradiction.
\end{pf}

\begin{lemma} \label{lem:4-cycle}
$G-y$ has a cycle of length four containing $x$.
\end{lemma}

\begin{pf}
Suppose that $x$ is not contained in any 4-cycle in $G-y$.
Then $d_{N(x)}(v)\le 1 \text{ for every } v\in V(G)-\{x,y\}$.

Let $G'$ be the graph obtained from $G$ by contracting $N[x]$ into a new vertex $x'$.
It is clear that $G'$ is connected and bipartite, and the minimum degree of $(G',x',y)$ is at least $k+1$ in $G'$.
If $G'$ is not 2-connected, then $x'$ is the unique cut-vertex of $G'$.
Let $H$ be the block of $G'$ containing $x'$ and $y$.
Note that $H=G'$ if $G'$ is 2-connected.

Suppose that $H$ is not an edge, then $(H,x',y)$ is 2-connected bipartite with minimum degree at least $k+1$.
By the induction hypothesis, $H$ contains $k$ paths $P_1',...,P_k'$ from $x'$ to $y$ with the length condition.
So $G-x$ contains $k$ paths $P_1,...,P_k$ from $N_G(x)$ to $y$ with the length condition.
Let $x_i$ be the end of $P_i$ contained in $N_G(x)$ for each $1 \leq i \leq k$.
By concatenating the edge $xx_i$ with $P_i$ for each $1 \leq i \leq k$, $G$ contains $k$ paths from $x$ to $y$ with the length condition, a contradiction.

Therefore, $H$ is an edge, which together with Lemma \ref{lem:xy_notadj} shows that $N_G(y) \subseteq N_G(x)$.
By \eqref{equ:dx<dy}, $N_G(x)=N_G(y)$.
We denote $N_G(x)$ by $N$.

Since $k \geq 3$ and $G$ is bipartite, $V(G)\neq N\cup \{x,y\}$.
So there exists a component $D$ of $G-N$ not containing $x$ and $y$.
Since $G$ is 2-connected, $|N_G(D)|\ge 2$.
Fixing a vertex $x''\in N_G(D)$,
let $G''$ be the graph obtained from $G[N_G[D]]$ by identifying $N_G(D)-x''$ into a new vertex $y''$.
Since $G$ is 2-connected and bipartite, $(G'',x'',y'')$ is also 2-connected and bipartite.
Since $d_{N}(v)\le 1 \text{ for every } v\in V(D)$, the minimum degree of $(G'',x'',y'')$ is at least $k+1$.
By induction, there exists a sequence of $k$ paths in $G''$ from $x''$ to $y''$ with the length condition.
So $G-\{x,y\}$ contains $k$ paths from $N$ to $N$ with the length condition.
By adding an edge between $x$ and $N$ and an edge between between $y$ and $N$ into each of these $k$ paths, we can obtain $k$ paths in $G$ from $x$ to $y$ with the length condition, a contradiction.
\end{pf}

\bigskip

The following notion is critical for the rest of the proof in this section.
Let $s$ be a positive integer.
A complete bipartite subgraph $Q$ of $G$ with bipartition $(Q_1,Q_2)$ is called an {\it $s$-core} if $x\in Q_2$, $y\notin V(Q)$, $|Q_1|\ge |Q_2|= s+1$,
and for every $v\in V(G)-(V(Q) \cup \{y\})$,
\begin{align}
\label{equ:d(v,Q)}
d_{Q_1}(v)\le s+1 \text{~~~ and ~~~} d_{Q_2}(v)\le s.
\end{align}
Since $G$ is bipartite, every vertex $v\in V(G)-(V(Q) \cup \{y\})$ is adjacent to at most one of $Q_1$ and $Q_2$, so $d_Q(v)=\max\{d_{Q_1}(v),d_{Q_2}(v)\}\le s+1$.

\medskip

The next lemma is straightforward but will be frequently used.
We omit the proof.

\begin{lemma} \label{lem:s-core-paths}
If $Q$ is an $s$-core in $G$, then for every $u\in Q_1$ there exist $s+1$ paths in $Q$ from $x$ to $u$ with lengths $1,3,\ldots,2s+1$, respectively, and for every $v\in Q_2-x$ there exist $s$ paths in $Q$ from $x$ to $v$ with lengths $2,4,\ldots,2s$, respectively.
\end{lemma}

\begin{lemma} \label{lem:s-core->d(v,Q1)}
$G$ contains an $s$-core $Q$ for some integer $s\ge 1$ such that the following hold.
Let $C$ be the component of $G-Q$ containing $y$.
If $G$ has an edge between $C$ and $Q_2-x$,
then for every $v\in V(G)-V(Q\cup C)$, $d_{Q_1}(v)\le s$ and thus $d_{Q}(v)\le s$.
\end{lemma}

\begin{pf}
Recall that $y$ is not adjacent to $x$ by Lemma \ref{lem:xy_notadj}.
By Lemma \ref{lem:4-cycle} there exists a 4-cycle in $G-y$ containing $x$.
Thus there exists a complete bipartite subgraph $Q$ of $G-y$ with bipartition $(Q_1,Q_2)$ such that $x \in Q_2$ and $\lvert Q_1 \rvert \geq \lvert Q_2 \rvert \geq 2$.
Let $C$ be the component of $G-V(Q)$ containing $y$.
We further choose $Q$ such that
	\begin{enumerate}
		\item[(a).] $\lvert Q_2 \rvert$ is maximum,
		\item[(b).] subject to (a), $Q_1$ is maximal, and
		\item[(c).] subject to (a) and (b), $\lvert V(C) \rvert$ is maximum.
	\end{enumerate}

Let $s=|Q_2|-1$.
We first prove that $Q$ is an $s$-core, which suffices to show \eqref{equ:d(v,Q)}.
Suppose to the contrary that there exists a vertex $v \in V(G)-(V(Q) \cup \{y\})$ satisfying that $d_{Q_1}(v)\ge s+2$ or $d_{Q_2}(v)\ge s+1$.
If $d_{Q_1}(v)\ge s+2$, then $|N_G(v)\cap Q_1|\ge s+2= |Q_2\cup\{v\}|$, and $G[(N_G(v)\cap Q_1)\cup Q_2\cup \{v\}]$ is a complete bipartite subgraph in $G-y$ with bipartition $(N_G(v)\cap Q_1,Q_2\cup \{v\})$, contradicting (a).
So $d_{Q_2}(v)\ge s+1$, that is $Q_2 \subseteq N_G(v)$.
Hence $(Q_1 \cup \{v\}, Q_2)$ is a complete bipartite subgraph of $G-y$, contradicting (b).
Therefore $Q$ is indeed an $s$-core.

Suppose that the lemma does not hold.
So by \eqref{equ:d(v,Q)}, there exists a vertex $v\in V(G)-V(Q\cup C)$ such that $\lvert N_G(v)\cap Q_1 \rvert = s+1$.
Assume that some vertex in $C$ is adjacent to a vertex $z\in Q_2-x$.
Let $Q_2' = Q_2 \cup \{v\}-\{z\}$, $Q_1' = \{a \in V(G): Q_2' \subseteq N_G(a)\}$, and $Q'=G[Q_1' \cup Q_2']$.
Since $y$ is not adjacent to $x$ in $G$, $y \not \in Q_1'$ and thus $y\notin V(Q')$.
Furthermore, $N_G(v)\cap Q_1\subseteq Q_1'$, so $Q'$ is a complete bipartite subgraph of $G-y$ containing $x$ with $|Q_1'|\ge s+1=|Q_2'|$, which also satisfies (a) and (b).
However, since $v$ is in a component of $G-V(Q)$ different from $C$, the component of $G-V(Q')$ containing $y$ contains $C$ and $z$.
This contradicts the choice of $Q$ as it violates (c).
This proves the lemma.
\end{pf}

\bigskip

In the rest of this section, $Q$ denotes the $s$-core mentioned in Lemma \ref{lem:s-core->d(v,Q1)}, and we let $C$ be the component of $G-V(Q)$ containing $y$.

Next we study the situation when there is an edge between $C$ and $Q_2-x$.
We will constantly use the following easy fact in the proofs:
if $A$ and $B$ are two arithmetic progressions with common difference two, then the elements of the set $\{a+b: a\in A, b\in B\}$ form an arithmetic progression of length $\lvert A \rvert + \lvert B \rvert -1$ with common difference two.

\begin{lemma} \label{lem:C-Q2->bad path}
If $C$ is adjacent in $G$ to some vertex $a\in Q_2-x$, then the following hold.
	\begin{enumerate}
		\item $G-V(C)$ does not contain $k$ paths from $x$ to $a$ satisfying the length condition.
		\item $G-V(C)$ does not contain $k-s+1$ paths from $Q_1$ to $Q_1$ internally disjoint from $V(Q)$ and satisfying the length condition.
		\item $G-V(C)$ does not contain $k-s+2$ paths from $Q_1$ to $Q_2-\{x,a\}$ internally disjoint from $V(Q)$ and satisfying the length condition.
		\item $G-V(C)$ does not contain $k-s+1$ paths from $Q_1$ to $\{x,a\}$ internally disjoint from $V(Q)$ and satisfying the length condition.
	\end{enumerate}
\end{lemma}

\begin{pf}
Suppose that $G-V(C)$ contains $k$ paths from $x$ to $a$ satisfying the length condition.
Then concatenating each path with a fixed path in $G[V(C)\cup \{a\}]$ from $a$ to $y$,
we obtain $k$ paths in $G$ from $x$ to $y$ satisfying the length condition, a contradiction.



Suppose that $G-V(C)$ contains $k-s+1$ paths $P_1,P_2,\ldots,P_{k-s+1}$ from $Q_1$ to $Q_1$ internally disjoint from $V(Q)$ and satisfying the length condition.
For each $i$, let $u_i,v_i\in Q_1$ be the two ends of $P_i$.
Then $Q-\{v_i,a\}$ contains $s$ paths from $x$ to $u_i$ with length $1,3,\ldots, 2s-1$, respectively.
By concatenating these $s$ paths with $P_i$ and the edge $v_ia$ for all $1\leq i\leq k-s+1$,
we obtain $k$ paths in $G-V(C)$ from $x$ to $a$ with the length condition,
a contradiction.

Suppose that $G-V(C)$ contains $k-s+2$ paths $P_1,P_2,\ldots,P_{k-s+2}$ from $Q_1$ to $Q_2-\{x,a\}$ internally disjoint from $V(Q)$ and satisfying the length condition.
For each $i$, let $u_i\in Q_2-\{x,a\}$ and $v_i\in Q_1$ be the ends of $P_i$. Then $Q-\{v_i,a\}$ contains $s-1$ paths from $x$ to $u_i$ with length $2,4,\ldots,2s-2$, respectively.
By concatenating these $s-1$ paths with $P_i$ and the edge $v_ia$ for all $1\leq i\leq k-s+2$, we obtain $k$ paths in $G-V(C)$ from $x$ to $a$ with the length condition, a contradiction.

Suppose that $G-V(C)$ contains $k-s+1$ paths $P_1,P_2,\ldots,P_{k-s+1}$
from $Q_1$ to $\{x,a\}$ internally disjoint from $V(Q)$ and satisfying the length condition.
For each $i$, let $u_i\in Q_1$ and $v_i\in \{x,a\}$ be the ends of $P_i$.
Then $Q-v_i$ contains $s$ paths from $u_i$ to $\{x,a\}-v_i$ with lengths $1,3,\ldots,2s-1$, respectively.
Concatenating these $s$ paths with $P_i$ for all $1\leq i\leq k-s+1$,
this gives rise to $k$ paths in $G-V(C)$ from $x$ to $a$ with the length condition,
a contradiction.
\end{pf}

\begin{lemma}\label{lem:C-Q2->N(Q1)}
If $C$ is adjacent in $G$ to some vertex $a\in Q_2-x$,
then 
$N_G(Q_1)\subseteq Q_2\cup V(C)$.
\end{lemma}

\begin{pf}
%
%
Suppose that $N_G(Q_1)\not \subseteq Q_2\cup V(C)$.
Then there is a component $D$ of $G-V(Q)$ other than $C$ with $\lvert N_G(D)\cap Q_1 \rvert \ge 1$.
Since $Q$ contains $s$ paths from $x$ to $a$ with the length condition, $s \leq k-1$ by Lemma \ref{lem:C-Q2->bad path}.

\medskip

{\bf Claim 1:} If $B$ is an end-block of $D$, then $N_G(B-b)\cap (Q_1\cup \{x,a\})\neq \emptyset$, where $b$ is the cut-vertex of $D$ contained in $B$.

{\bf Proof of Claim 1.}
Suppose to the contrary that $N_G(B-b)\cap V(Q)\subseteq Q_2-\{x,a\}$.
Since $G$ is 2-connected, we have $|V(Q_2)-\{x,a\}|\ge 1$ and thus $s\ge 2$.
Let $G_1$ be the graph obtained from $G[V(B) \cup (N_G(B-b)\cap V(Q))]$ by identifying $N_G(B-b)\cap V(Q)$ into a vertex $x_1$.
So $(G_1,x_1,b)$ is 2-connected bipartite and has minimum degree at least $(k+1)-(s-2)$.
By induction $G_1$ has $k-s+2$ paths from $x_1$ to $b$ with the length condition.
There is a path in $N_G[D-V(B-b)]$ from $b$ to $Q_1$.
So $G-V(C)$ has $k-s+2$ paths from $Q_2-\{x,a\}$ to $Q_1$ internally disjoint from $V(Q)$ and satisfying the length condition, contradicting Lemma \ref{lem:C-Q2->bad path}.~
$\Box$

\medskip

{\bf Claim 2:} $N_G(D)\cap \{x,a\}=\emptyset$.

{\bf Proof of Claim 2.}
Suppose that $N_G(D)\cap \{x,a\} \neq \emptyset$.
Let $G_2$ be the graph obtained from $N_G[D]-(Q_2-\{x,a\})$ by identifying $N_G(D)\cap \{x,a\}$ into a vertex $x_2$ and identifying $N_G(D)\cap Q_1$ into a vertex $y_2$.
For every $v\in V(G_2)-\{x_2,y_2\}$, $d_Q(v)\le s$ by Lemma \ref{lem:s-core->d(v,Q1)}, and $v$ is adjacent to at most one of $Q_1$ and $Q_2$.
If $v$ is not adjacent to $Q_1$ or $Q_2$, then $d_{G_2}(v) \geq k+1$; if $v$ is adjacent to $Q_1$, it is clear that $d_{G_2}(v)\ge (k+1)-(s-1)$;
if $v$ is adjacent to $Q_2$ but not to any one of $x,a$, then $d_Q(v)\le s-1$, implying that $d_{G_2}(v)\ge (k+1)-(s-1)$; otherwise $v$ is adjacent to at least one of $x,a$, then $d_{G_2}(v)\ge (k+1)-(s-1)$.
Therefore $(G_2,x_2,y_2)$ has minimum degree at least $k-s+2$.
By Claim 1, every end-block of $G_2$ contains at least one of $x_2,y_2$ as a non-cut-vertex,
so $(G_2,x_2,y_2)$ is 2-connected and bipartite.
By induction, $G_2$ contains $k-s+1$ paths from $x_2$ to $y_2$ satisfying the length condition.
So $G-V(C)$ contains $k-s+1$ paths from $\{x,a\}$ to $Q_1$ internally disjoint from $V(Q)$ and satisfying the length condition, contradicting Lemma \ref{lem:C-Q2->bad path}.
$\Box$

\medskip

{\bf Claim 3:} $\lvert N_G(D)\cap Q_1 \rvert \ge 2$.

{\bf Proof of Claim 3.}
Suppose to the contrary that $\lvert N_G(D)\cap Q_1 \rvert \leq 1$.
By the choice of the component $D$, $N_G(D)\cap Q_1=\{x_3\}$ for some vertex $x_3$.
Since $G$ is 2-connected, Claim 2 implies that $|N_G(D)\cap (Q_2-\{x,a\})|\ge 1$, so $s\ge 2$.
Let $G_3$ be the graph obtained from $N_G[D]$ by identifying $N_G(D)\cap (Q_2-\{x,a\})$ into a vertex $y_3$.
In view of Claim 2, every end-block of $G_3$ contains at least one of $x_3,y_3$ as a non-cut-vertex,
so $(G_3,x_3,y_3)$ is 2-connected and bipartite.
For any $v\in V(G_3)-\{x_3,y_3\}$, if $v$ is adjacent to $Q_1$, then $d_{G_3}(v)=d_G(v)\ge k-s+3$;
otherwise $N_G(v)\cap Q\subseteq Q_2-\{x,a\}$, also implying $d_{G_3}(v)\ge (k+1)-(s-2)=k-s+3$.
By induction, $G_3$ contains $k-s+2$ paths from $x_3$ to $y_3$ with the length condition.
Hence, $G-V(C)$ contains $k-s+2$ paths from $Q_1$ to $Q_2-\{x,a\}$ internally disjoint from $V(Q)$ and satisfying the length condition, contradicting Lemma \ref{lem:C-Q2->bad path}.
$\Box$

\medskip

Fix a vertex $x_4\in N_G(D)\cap Q_1$.
Claim 3 ensures that $N_G(D)\cap Q_1-x_4 \neq \emptyset$.
Let $G_4$ be the graph obtained from $G[N_G[D]-Q_2]$ by identifying $N_G(D)\cap Q_1-x_4$ into a vertex $y_4$.
Recall Lemma \ref{lem:s-core->d(v,Q1)} that $d_Q(v)\le s$ for every $v\in V(D)$.
For every $v\in V(G_4)-\{x_4,y_4\}$ adjacent in $G$ to $Q$, if $v$ is adjacent to $Q_1$, then $d_{G_4}(v)\ge (k+1)-(s-1)$;
otherwise $v$ is adjacent to $Q_2$, so $d_{G_4}(v)\ge (k+1)-(s-1)$ by Claim 2.
Hence $(G_4,x_4,y_4)$ has minimum degree at least $k-s+2$.
By Claims 1 and 2, every end-block of $G_4$ contains at least one of $x_4,y_4$ as a non-cut-vertex, so $(G_4,x_4,y_4)$ is 2-connected and bipartite.
By induction, $G_4$ contains $k-s+1$ paths from $x_4$ to $y_4$ satisfying the length condition.
So $G-V(C)$ contains $k-s+1$ paths from $Q_1$ to $Q_1$ internally disjoint from $V(Q)$ and satisfying the length condition, contradicting Lemma \ref{lem:C-Q2->bad path}.
\end{pf}

\begin{lemma} \label{lem:C>=2}
$C$ contains at least two vertices, and no vertex of $C-y$ is a leaf in $C$.
\end{lemma}

\begin{pf}
We first prove that no vertex of $C-y$ is a leaf in $C$.
Suppose that $C$ has a leaf $z\in V(C-y)$.
If $z$ is adjacent to $Q_1$, by \eqref{equ:d(v,Q)} we have $s+1\ge d_{Q_1}(z)\ge k$,
so by Lemma \ref{lem:s-core-paths}, there are $k$ paths in $V(Q)$ from $x$ to $N_G(z)\cap Q_1$ with lengths $1,3,\ldots,2k-1$, respectively,
which can be easily extended to $k$ paths in $G$ from $x$ to $y$ with the length condition.
Hence $z$ is adjacent to $Q_2$. By Lemma \ref{basis} and \eqref{equ:d(v,Q)}, we have $s\geq d_{Q_2}(z)\geq k\geq 3$,
so there is a vertex $a\in N_G(z)\cap Q_2-x$.
By Lemma \ref{lem:s-core-paths}, there are $k$ paths in $Q$ from $x$ to $a$ with the length condition,
contradicting Lemma \ref{lem:C-Q2->bad path}.

It suffices to show that $C$ has at least two vertices.
We suppose for a contradiction that $C$ consists of one vertex, i.e., $V(C) = \{y\}$.

\medskip

{\bf Claim 1:} $N_G(x)=N_G(y)=Q_1$ and $V(G) \neq V(Q\cup C)$.

{\bf Proof of Claim 1:}
If $y$ is adjacent in $G$ to a vertex $a \in Q_2-x$, then $N_G(Q_1)\subseteq Q_2\cup \{y\}$ by Lemma \ref{lem:C-Q2->N(Q1)}.
Since $G$ is bipartite, $N_G(Q_1) \subseteq Q_2$, so $s \geq k$.
Then by Lemma \ref{lem:s-core-paths}, $Q$ contains $k$ paths from $x$ to $a$ with the length condition, contradicting Lemma \ref{lem:C-Q2->bad path}.
Hence $N_G(y) \subseteq Q_1 \cup \{x\}$.
But $x$ is not adjacent to $y$, so $N_G(y) \subseteq Q_1 \subseteq N_G(x)$.
By the assumption \eqref{equ:dx<dy}, the degree of $x$ in $G$ is at most the degree of $y$ in $G$.
This proves that $N_G(x)=N_G(y)=Q_1$.

Similarly, if $V(G)=V(Q\cup C)$, then $N_G(Q_1)=Q_2\cup \{y\}$ and $s \geq k-1$.
Let $z \in N_G(y) \cap Q_1$.
By Lemma \ref{lem:s-core-paths}, $Q$ contains $s+1\geq k$ paths from $x$ to $z$ with the length condition, a contradiction.
Therefore, $V(G) \neq V(Q\cup C)$.
$\Box$

%
%

\medskip

{\bf Claim 2:} $\lvert Q_1 \rvert \ge 3$.

{\bf Proof of Claim 2:}
Suppose $\lvert Q_1 \rvert \leq 2$, then $\lvert Q_1 \rvert=2$ and $s=1$.
Let $Q_1=\{u,w\}$, $Q_2=\{v,x\}$ and $G_1=G-\{x,y\}$.
Note that $G_1$ is connected.
By Claim 1, $N_G(x)=N_G(y)=\{u,w\}$, so $(G_1,u,w)$ has minimum degree at least $k+1$ in $G_1$.
If $G_1$ is 2-connected, then $(G_1,u,w)$ is 2-connected.
Otherwise, since $G$ is 2-connected and $N_G(x)=N_G(y)=\{u,w\}$, $u$ and $w$ are in different end-blocks of $G_1$;
since $u,w\in N_G(v)$, $v$ is the cut-vertex of $G_1$ contained in both of the two end-blocks of $G_1$.
So $(G_1,u,w)$ is 2-connected in either case.
Therefore, by induction, $G_1$ has $k$ paths from $u$ to $w$ satisfying the length condition.
Concatenating them with $xu,wy$ gives $k$ path in $G$ from $x$ to $y$ satisfying the length condition.
$\Box$

\medskip

Let $u\in Q_1$ and $v\in Q_2-x$ be fixed. Then any vertex in $G-\{u,v\}$ other than $x,y$ has degree at least $k$ in $G-\{u,v\}$.
If $G-\{u,v\}$ is 2-connected, then $G-\{u,v\}$ contains $k-1$ paths from $x$ to $y$ with the length condition.
Among these paths, let $R$ be the longest one such that $w\in Q_1-u$ is the end of $R-y$ other than $x$.
These $k-1$ paths from $x$ to $y$ together with the path $(R-y)\cup wvuy$ are $k$ paths in $G$ from $x$ to $y$ with the length condition.
Hence $G-\{u,v\}$ is not 2-connected and contains at least two end-blocks.

Suppose that there exists a component $H$ of $G-\{u,v\}$ disjoint from $Q \cup C$.
Since $G$ is 2-connected, $(G[V(H) \cup \{u,v\}],u,v)$ is 2-connected bipartite and has minimum degree at least $k+1$. 
So $G$ contains $k$ paths from $u$ to $v$ internally disjoint from $V(Q) \cup \{y\}$ with the length condition.
By concatenating $xu$ and $vwy$ with each path, where $w$ is a vertex in $Q_1-u$, we obtain $k$ paths in $G$ from $x$ to $y$ with the length condition, a contradiction.
Therefore, $G-\{u,v\}$ is connected.

By Claims 1 and 2, $G[V(Q\cup C)]-\{u,v\}$ is 2-connected.
So there is an end-block $B$ of $G-\{u,v\}$ with the cut-vertex $b$ such that $(B-b)\cap ((Q\cup C)-\{u,v\})=\emptyset$.
Since $G-\{u,v\}$ is connected, there exists a path $P$ in $G-\{u,v\}$ from $b$ to
some vertex $z\in V(Q \cup C)-\{u,v\}$ internally disjoint from $(B \cup Q \cup C)-\{u,v\}$.
Note that $z\notin \{x,y\}$ as $N_G(x)=N_G(y)=Q_1$. So $z\in V(Q)-\{u,v,x\}$.

Note that $N_G(B-b) \subseteq \{u,v,b\}$.
Suppose that $u \not \in N_G(B-b)$.
Then $(G[V(B) \cup \{v\}],v,b)$ is 2-connected bipartite with minimum degree at least $k+1$.
So induction ensures that $G[V(B)\cup \{v\}]$ contains $k$ paths $P_1,P_2,...,P_k$ from $v$ to $b$ with the length condition.
If $z\in Q_1$, let $P'=P\cup \{zy\}$; if $z\in Q_2$, fix a vertex $w\in Q_1-u$ and let $P'=P\cup \{zw,wy\}$.
So in either case $P'$ is a path from $b$ to $y$ and internally disjoint from $B\cup \{u,v,x\}$.
By concatenating $P_i$ with $xuv$ and $P'$ for each $1\leq i\leq k$,
we obtain $k$ paths in $G$ from $x$ to $y$ with the length condition.
Therefore, $u\in N_G(B-b)$ and hence $(G[V(B) \cup \{u\}],u,b)$ is 2-connected.

Since $(G[V(B) \cup \{u\}],u,b)$ is 2-connected bipartite with minimum degree at least $k$, $G[V(B) \cup \{u\}]$ contains $k-1$ paths from $u$ to $b$ with the length condition.
By concatenating these paths with $P$, this gives a sequence of $k-1$ paths $R_1,R_2,...,R_{k-1}$ in $G-v$ from $u$ to $z$ internally disjoint from $V(Q \cup C)$ with the length condition.
If $z \in Q_1$, then by Claim 2, there exists a vertex $w \in Q_1-\{u,z\}$, and we let $R_k$ be the path obtained from $R_{k-1}$ by concatenating $zvw$.
Then $R_1, R_2,...,R_k$ form a sequence of $k$ paths in $G-\{x,y\}$ from $Q_1$ to $Q_1$ with the length condition, which, by Claim 1, can be easily extended to $k$ path in $G$ from $x$ to $y$ with the length condition.
Thus $z\in Q_2$. By Claim 2, there exist two distinct vertices $w,w'\in Q_1-u$.
For each $1 \leq i \leq k-1$, let $R_i'$ be the path obtained from $R_i$ by concatenating $xu$ and $zwy$; and let $R_k'$ be the path obtained from $R_{k-1}$ by concatenating $xu$ and $zwvw'y$.
Therefore, $R_1',R_2',..., R_k'$ form a sequence of $k$ paths in $G$ from $x$ to $y$ with the length condition.
This proves the lemma.
\end{pf}

\begin{lemma} \label{middle to not y exists}
$G$ has an edge between $Q_1$ and $C-y$.
\end{lemma}

\begin{pf}
Note that $C-y \neq \emptyset$ by Lemma \ref{lem:C>=2}.
Suppose to the contrary that $N_G(C-y)\cap Q_1=\emptyset$.

We claim that $N_G(C-y)\cap (Q_2-x)=\emptyset$.
Otherwise, $C-y$ is adjacent to some vertex $a\in Q_2-x$.
By Lemma \ref{lem:C-Q2->N(Q1)} and the assumption $N_G(C-y)\cap Q_1=\emptyset$,
it follows that $N_G(Q_1)\subseteq Q_2\cup \{y\}$.
So for some $u\in Q_1$, $N_G(u)\subseteq Q_2\cup \{y\}$.
This implies that $s\ge k-1$, and if $s=k-1$, then $uy\in E(G)$.
If $s\ge k$, by Lemma \ref{lem:s-core-paths}, there are at least $k$ paths in $Q$ from $x$ to $a$ with the length condition,
contradicting Lemma \ref{lem:C-Q2->bad path}.
So $s=k-1$ and thus $uy\in E(G)$.
Again by Lemma \ref{lem:s-core-paths}, there are at least $k$ paths in $Q$ from $x$ to $u$ with the length condition.
Concatenating them with $uy$ gives $k$ paths from $x$ to $y$ with the length condition, a contradiction. This proves that $N_G(C-y)\cap (Q_2-x)=\emptyset$.

Therefore, $N_G(C-y)=\{x,y\}$.
Since $G$ is 2-connected, $(G[V(C) \cup \{x\}],x,y)$ is 2-connected bipartite and has minimum degree at least $k+1$.
By the induction hypothesis, $G$ contains $k$ paths from $x$ to $y$ satisfying the length condition.
\end{pf}

\begin{lemma} \label{lem:bad path from y to Q}
$G$ does not contain $k-s$ paths from $y$ to $Q_1$ internally disjoint from $V(Q)$ with the length condition nor $k-s+1$ paths from $y$ to $Q_2-x$ internally disjoint from $V(Q)$ with the length condition.
\end{lemma}

\begin{pf}
Suppose to the contrary that there exist $k-s$ paths $P_1,\dots, P_{k-s}$ in $G$ from $y$ to $Q_1$ internally disjoint from $V(Q)$ and satisfying the length condition.
For each $1\leq i\leq k-s$, let $u_i\in Q_1$ be the end of $P_i$ other than $y$.
By Lemma \ref{lem:s-core-paths}, $Q$ contains $s+1$ paths from $x$ to $u_i$ with lengths $1,3,...,2s+1$, respectively.
Then concatenating these $s+1$ paths with $P_i$ for each $1 \leq i \leq k-s$ leads to $k$ paths in $G$ from $x$ to $y$ with the length condition, a contradiction.

Suppose to the contrary that there exist $k-s+1$ paths $R_1,\dots, R_{k-s+1}$ in $G$ from $y$ to $Q_2-x$ internally disjoint from $V(Q)$ and satisfying the length condition.
For each $1\leq j\leq k-s+1$, let $v_j\in Q_2-x$ be the end of $R_j$ other than $y$.
By Lemma \ref{lem:s-core-paths}, $Q$ contains $s$ paths from $x$ to $v_j$ with lengths $2,4,\ldots,2s$, respectively.
Then concatenating these $s$ paths with $R_j$ for each $1\le j\le k-s+1$ leads to $k$ paths in $G$ from $x$ to $y$ with the length condition, a contradiction.
\end{pf}

\bigskip

We say that an end-block $B$ of $C$ is {\it feasible} if $y\notin V(B-b)$, where $b$ is the cut-vertex of $C$ contained in $B$.

\begin{lemma} \label{neighbor in Q2}
$s=1$, and $C$ is not 2-connected.
Moreover, if $B$ is a feasible end-block of $C$ with the cut-vertex $b$, then $B$ is 2-connected and $N_G(B-b)= Q_2 \cup \{b\}$.
\end{lemma}

\begin{pf}
Recall that $C$ contains at least two vertices, and no vertex of $C-y$ is a leaf in $C$ by Lemma \ref{lem:C>=2}.
So every feasible end-block of $C$ is 2-connected.

\medskip

{\bf Claim 1:} $C$ is not 2-connected, and for each feasible end-block $B$ of $C$ with cut-vertex $b$, $N_G(B-b)\cap Q_1=\emptyset$.

{\bf Proof of Claim 1.}
Suppose to the contrary.
So either $C$ is 2-connected,
or there is an end-block $B$ of $C$ with cut-vertex $b$ such that $y\notin V(B-b)$ and $B-b$ is adjacent in $G$ to $Q_1$.
In the former case, define $B'=C$ and $b'=y$, so $B'-b'$ is adjacent to $Q_1$ by Lemma \ref{middle to not y exists}; in the latter case, define $B'=B$ and $b'=b$, so $B'-b'$ is adjacent to $Q_1$ by the assumption.
Note that there is a path $P$ in $C$ from $b'$ to $y$ internally disjoint from $B'$.
Let $X= N_G(B'-b')\cap Q_1$ and define $G_1$ to be the graph obtained from $G[B'\cup X]$ by identifying $X$ into a vertex $x_1$.
By \eqref{equ:d(v,Q)}, $(G_1,x_1,b')$ has minimum degree at least $k+1-s$.
Since $(G_1,x_1,b')$ is 2-connected and bipartite, by induction $G_1$ has $k-s$ paths from $b'$ to $x_1$ with the length condition.
By concatenating with the path $P$, it is easy to obtain $k-s$ paths in $G$ from $y$ to $Q_1$ internally disjoint from $V(Q)$ and satisfying the length condition, contradicting Lemma \ref{lem:bad path from y to Q}.
$\Box$

\medskip

Claim 1 implies that feasible end-blocks of $C$ exist.
Let $B$ be an arbitrary feasible end-block of $C$, and let $b$ be the cut-vertex of $C$ contained in $B$.

\medskip

{\bf Claim 2:} $N_G(B-b)\cap (Q_2-x)\neq \emptyset$.

{\bf Proof of Claim 2.}
Suppose to the contrary that $N_G(B-b)\cap (Q_2-x)= \emptyset$.
By Claim 1 and the 2-connectivity of $G$, $N_G(B-b)=\{b,x\}$.
Define $G_2=G[V(B)\cup \{x\}]$.
Since $(G_2,x,b)$ is 2-connected bipartite and has minimum degree at least $k+1$,
$G_2$ has $k$ paths from $x$ to $b$ with the length condition.
By concatenating them with a fixed path in $C-V(B-b)$ from $b$ to $y$, we obtain $k$ paths in $G$ from $x$ to $y$ with the length condition.
$\Box$

\medskip

Finally, we shall prove that $s=1$ and $N_G(B-b) =Q_2 \cup \{b\}$.
Suppose that either $s\ge 2$, or $s=1$ but $N_G(B-b) \neq Q_2 \cup \{b\}$.
Note that the latter case implies that $N_G(B-b)=\{b\} \cup (Q_2-x)$ by Claims 1 and 2.
Define $G_3$ to be the graph obtained from $G[V(B)\cup (Q_2-x)]$ by identifying $Q_2-x$ into vertex $a'$.
Claim 2 implies that $(G_3,a',b)$ is 2-connected and bipartite.

We show that every vertex $v\in V(G_3)-\{a',b\}$ has degree at least $k-s+2$ in $G_3$.
Note that $v$ has at most $s$ neighbors in $Q_2$ by \eqref{equ:d(v,Q)} and no neighbor in $Q_1$ by Claim 1.
If $s\ge 2$ and $v$ has at most $s-1$ neighbors in $Q_2$, then it is clear that $d_{G_3}(v)\ge (k+1)-(s-1)$.
If $s\ge 2$ and $v$ has exactly $s$ neighbors in $Q_2$, then at least one of them is in $Q_2-x$ and thus $d_{G_3}(v)\ge (k+1)-(s-1)$.
It remains to consider $s=1$. In this case, as $x\notin N_G(B-b)$, it is easy to see that $d_{G_3}(v)\ge k+1$.
Therefore, $(G_3,a',b)$ has minimum degree at least $k-s+2$.

By induction, $G_3$ has $k-s+1$ paths from $a'$ to $b$ with the length condition.
Concatenating them with a fixed path in $C-V(B-b)$ from $b$ to $y$, we can obtain $k-s+1$ paths in $G$ from $y$ to $Q_2-x$
internally disjoint from $V(Q)$ and satisfying the length condition, contradicting Lemma \ref{lem:bad path from y to Q}.
\end{pf}

\bigskip

By Lemma \ref{neighbor in Q2}, $C$ has at least two end-blocks, but at most one of them contains $y$ as a non-cut-vertex. 
So there is at least one feasible end-block of $C$.
We also see that $Q_2-x$ contains exactly one vertex from Lemma \ref{neighbor in Q2}.
In the rest of this section, we denote this vertex by $a$.
Namely, $Q_2=\{a,x\}$.

\begin{lemma} \label{feasible k-1}
Let $B$ be a feasible end-block of $C$ with the cut-vertex $b$.
For each vertex $u$ in $Q_2=\{a,x\}$, $G[V(B)\cup \{u\}]$ has $k-1$ paths from $u$ to $b$ with the length condition.
\end{lemma}

\begin{pf}
Define $G'=G[V(B)\cup \{u\}]$.
So $(G',u,b)$ is 2-connected and bipartite.
By Lemma \ref{neighbor in Q2}, $(G',u,b)$ has minimum degree at least $k$.
By induction, $G'$ has $k-1$ paths from $u$ to $b$ with the length condition.
\end{pf}

\bigskip

We complete the proof of Theorem \ref{main bipartite} in the coming last lemma of this section.

\begin{lemma}\label{lem:final step}
$G$ is not a counterexample of Theorem \ref{main bipartite}.
\end{lemma}

\begin{pf}
Define $N=N_G(Q_1) \cap V(C-y)$. Lemma \ref{middle to not y exists} implies that $N\neq \emptyset$.
Let $B_1, B_2,..., B_t$ be all feasible end-blocks of $C$, and let $b_i$ be the cut-vertex of $C$ contained in $B_i$ for each $i$.
Let $C'$ be obtained from $C$ by deleting $V(B_i-b_i)$ for all $i$.
By Lemma \ref{neighbor in Q2} and the definition of feasible end-blocks, $C'$ is connected and contains $N\cup \{y,b_1,b_2,...,b_t\}$.

\medskip

{\bf Claim 1:} There exists $c \in V(C')$ such that no path in $C'-c$ is from $N\cup \{y\}$ to $\{b_1,b_2,...,b_t\}$.

{\bf Proof of Claim 1.}
Suppose to the contrary that there exist two disjoint paths $P_1,P_2$ in $C'$ from $N\cup \{y\}$ to $\{b_1,b_2,...,b_t\}$.
Since $C'$ is connected, we may assume that $y$ is an end of one of $P_1,P_2$, say $P_1$, by rerouting paths.
Denote the end of $P_2$ in $N$ by $w$.
By symmetry, we may without loss of generality assume that the ends of $P_1,P_2$ in $\{b_1,b_2,...,b_t\}$ are $b_1$ and $b_2$, respectively. 
By Lemma \ref{feasible k-1}, there exist a sequence of $k-1$ paths $R_1,R_2,...,R_{k-1}$ in $G[V(B_1) \cup \{a\}]$ from $a$ to $b_1$ with the length condition and a sequence of $k-1$ paths $L_1,L_2,...,L_{k-1}$ in $G[V(B_2) \cup \{x\}]$ from $x$ to $b_2$ with the length condition.
Let $w' \in Q_1\cap N_G(w)$.
Since $k\ge 3$ by Lemma \ref{basis}, for all $i,j \in \{1,2,...,k-1\}$, the paths $L_i\cup P_2\cup ww'a\cup R_j \cup P_1$ give rise to at least $2k-3 \geq k$ paths in $G$ from $x$ to $y$ satisfying the length condition, a contradiction.
$\Box$

\medskip

{\bf Claim 2:} There exists an end-block $B_y$ of $C$ with cut-vertex $b_y$ such that $y\in V(B_y-b_y)$.

{\bf Proof of Claim 2.} Otherwise, all end-blocks of $C$ are feasible.
By Claim 1, there exist a cut-vertex $c$ of $C'$ and two subgraphs $C_1,C_2$ of $C'$ such that $C'=C_1\cup C_2$ and $V(C_1)\cap V(C_2)=\{c\}$, where $N\cup \{y\}\subseteq C_1$ and $\{b_1,b_2,...,b_t\}\subseteq C_2$.
But $C_2$ contains all cut-vertices of $C$ contained in some end-blocks of $C$, a contradiction.
$\Box$

\medskip

{\bf Claim 3:} For every $v\in V(C-y)$, either $d_Q(v)\leq 1$ or $v$ is a cut-vertex of $C$ separating $y$ and all feasible end-blocks of $C$.

{\bf Proof of Claim 3.}
Suppose to the contrary that there exist a vertex $v\in V(C-y)$ with $d_Q(v)\geq 2$ and a feasible end-block $B$ of $C$ with cut-vertex $b$ such that $C-v$ has a path $L$ from $y$ to $b$ internally disjoint from $B$. Since $s=1$ by Lemma \ref{neighbor in Q2},
\eqref{equ:d(v,Q)} ensures that $v$ is adjacent to two distinct vertices in $Q_1$, say $u_1,u_2$.
By Lemma \ref{feasible k-1}, there exists a sequence of $k-1$ paths $P_1,P_2,...,P_{k-1}$ in $G[V(B)\cup \{a\}]$ from $a$ to $b$ with the length condition.
Then $xu_1a \cup P_i \cup L$ for all $1\leq i \leq k-1$ together with $xu_2vu_1a \cup P_{k-1} \cup L$ are $k$ paths in $G$ from $x$ to $y$ with the length condition, a contradiction.
$\Box$

\medskip

Fix a feasible end-block $B$ of $C$, and let $b$ be the cut-vertex of $C$ contained in $B$.
By Lemma \ref{feasible k-1}, there exists a sequence of $k-1$ paths $P_1,P_2,..., P_{k-1}$ in $G[V(B) \cup \{x\}]$ from $x$ to $b$ with the length condition.
Concatenating them with a fixed path in $C$ from $b$ to $b_y$, we obtain a sequence of $k-1$ paths $R_1,R_2,...,R_{k-1}$ in $G[(V(C) \cup \{x\}) -V(B_y-b_y)]$ from $x$ to $b_y$ with the length condition.

\medskip

{\bf Claim 4:} $B_y$ is an edge $yb_y$.

{\bf Proof of Claim 4.}
Suppose to the contrary that $B_y$ is 2-connected.
For every $v\in V(B_y)-\{y,b_y\}$, $v$ is not a cut-vertex of $C$ separating $y$ and feasible end-blocks of $C$, so $d_Q(v)\le 1$ by Claim 3.
So $(B_y,y,b_y)$ is 2-connected bipartite with minimum degree at least $k$.
By induction, $B_y$ contains $k-1$ paths from $y$ to $b_y$ with the length condition.
Concatenating these $k-1$ paths with $R_i$ for each $1\leq i\leq k-1$,
we obtain $2k-3 \geq k$ paths in $G$ from $x$ to $y$ with the length condition, a contradiction.
$\Box$

\medskip

Suppose that $b_y$ is adjacent in $G$ to a vertex $z\in Q_1$.
If $y$ is adjacent to $Q_1$, 
then Claim 4 will force an odd cycle in $G$, a contradiction as $G$ is bipartite.
So $N_G(y)\subseteq Q_2\cup \{b_y\}$.
Since $G$ is 2-connected and $xy\notin E(G)$, $N_G(y)=\{a,b_y\}$.
Then $R_i \cup b_yy$ for all $1\leq i \leq k-1$ together with $R_{k-1} \cup b_yzay$ form $k$ paths in $G$ from $x$ to $y$ with the length condition, a contradiction.
Therefore, $b_y$ is not adjacent to $Q_1$, that is, $b_y\notin N$.
Also by \eqref{equ:d(v,Q)}, $d_Q(b_y)\le 1$.

Let $W$ be a block of $C-y$ containing $b_y$. Since $d_Q(b_y)\le 1$, we have $d_{C-y}(b_y)\geq k-1\ge 2$, so $W$ is 2-connected.
If $W=B_i$ for some $i$, then $V(C)=V(B_i)\cup \{y\}$, $b_y=b_i$, and $b_y$ is adjacent to $Q_1$ by Lemmas \ref{middle to not y exists} and \ref{neighbor in Q2}, a contradiction.
So $W$ is not an end-block of $C$ and thus $W\cup \{y\}\subseteq C'$.

Since $W$ is 2-connected and $b_y\notin N$, Claim 1 implies that there exists a cut-vertex of $C'$ separating $N\cup W\cup \{y\}$ and $\{b_1,b_2,...,b_t\}$.
Hence $C-y$ has a cut-vertex separating $W$ and all feasible end-blocks of $C$.
Note that every cut-vertex of $C-y$ contained in $W$ has a path to some feasible end-block of $C$ internally disjoint from $W$.
Therefore, $W$ has the unique cut-vertex $w$ of $C-y$.

For every $v\in V(W)-\{w,b_y\}$, since $v$ is not a cut-vertex of $C$ separating $y$ and all feasible end-blocks of $C$, we have $d_Q(v)\leq 1$ by Claim 3.
This together with $d_Q(b_y)\le 1$ imply that
$(G[V(W) \cup \{y\}],w,y)$ is 2-connected bipartite with minimum degree at least $k$.
By induction, there exists a sequence of $k-1$ paths $L_1,L_2,...,L_{k-1}$ in $G[V(W) \cup \{y\}]$ from $w$ to $y$ with the length condition.
Recall the $k-1$ paths $P_1,P_2,...,P_{k-1}$ in $G[V(B)\cup\{x\}]$ from $x$ to $b$.
Let $R$ be a path in $C$ from $b$ to $w$ internally disjoint from $B\cup W\cup \{y\}$.
Then for all $i,j \in \{1,2,..., k-1\}$, the paths $P_i\cup R\cup L_j$ give rise to $2k-3 \geq k$ paths in $G$ from $x$ to $y$ with the length condition, a contradiction.
\end{pf}

\medskip

This proves Theorem \ref{main bipartite}, which implies Theorem \ref{thm:bipartite_main}.

\section{Consecutive paths in general graphs} \label{sec:path_general}
\noindent The following two lemmas extend Theorem \ref{main bipartite} from bipartite graphs to general graphs, which will be extensively used in the coming sections for finding cycles.
\begin{lemma} \label{path general graph}
Let $(G,x,y)$ be a 2-connected rooted graph.
If the minimum degree of $(G,x,y)$ is at least $k+1$, then $G$ contains $\lfloor k/2\rfloor$ paths from $x$ to $y$ satisfying the length condition.
\end{lemma}

\begin{pf}
%
Let $G'$ be a spanning bipartite subgraph of $G$ with maximum number of edges.
So for every vertex $v\in V(G)$, we have $d_{G'}(v)\geq \lceil d_G(v)/2\rceil$.
Hence, every vertex of $G'$ other than $x,y$ has degree at least $\lfloor k/2\rfloor+1$.
By the maximality, $G'$ is connected.

Suppose that there exists an end-block $B$ of $G'$ such that $V(B-b)\cap \{x,y\}=\emptyset$, where $b$ is the cut-vertex of $G'$ contained in $B$.
There exists a path $P$ in $G-(B-b)$ from $b$ to $\{x,y\}$ as $G'$ is connected.
Since $(G,x,y)$ is 2-connected, 
there exist two disjoint paths in $G$ from $V(B)$ to $\{x,y\}$ internally disjoint from $V(B)$.
Rerouting these two paths by the path $P$, we can further obtain two disjoint paths $P_1,P_2$ in $G$ from $V(B)$ to $\{x,y\}$ internally disjoint from $V(B)$ such that $b$ is an end of $P_1$ or $P_2$, say $P_1$.
We denote the end of $P_2$ in $B$ by $u$.
Every vertex in $V(B-b)$ has degree at least $\lfloor k/2\rfloor+1$ in $B$, so $B$ is 2-connected bipartite with minimum degree at least $\lfloor k/2 \rfloor+1$.
By Theorem \ref{thm:bipartite_main}, $B$ contains $\lfloor k/2\rfloor$ paths from $b$ to $u$ with the length condition.
By concatenating each of them with the paths $P_1,P_2$, we obtain $\lfloor k/2\rfloor$ paths in $G$ from $x$ to $y$ satisfying the length condition.

Therefore, every end-block of $G'$ contains at least one of $x,y$ as a non-cut-vertex.
So $(G',x,y)$ is 2-connected bipartite with minimum degree at least $\lfloor k/2\rfloor+1$,
by Theorem \ref{main bipartite} there exist $\lfloor k/2\rfloor$ paths in $G'$ (and hence in $G$) from $x$ to $y$ satisfying the length condition.
\end{pf}

\begin{lemma}\label{lem:path with one exceptional vertex}
Let $G$ a 2-connected graph and $x,y,v$ be distinct vertices of $G$.
If every vertex of $G$ other than $v$ has degree at least $k+1$, then $G$ contains $\lfloor (k-1)/2\rfloor$ paths from $x$ to $y$ satisfying the length condition.
\end{lemma}

\begin{pf}
There is nothing to prove when $k\leq 2$, so we may assume that $k\geq 3$.
Note that $G-v$ is connected and has minimum degree at least $k$.
If $G-v$ is 2-connected, then it follows from Lemma \ref{path general graph}.
Hence we may assume that $G-v$ is not 2-connected.
Then any end-block of $G-v$ is 2-connected and has a non-cut-vertex adjacent to $v$ in $G$.

Let $B$ be an arbitrary end-block of $G-v$, and let $b$ be the cut-vertex of $G-v$ contained in $B$.
Suppose that $\lvert V(B-b)\cap \{x,y\} \rvert=1$.
Without loss of generality, we may assume that $x\in V(B-b)$.
By Lemma \ref{path general graph},
$B$ has $\lfloor (k-1)/2\rfloor$ paths from $x$ to $b$ with the length condition.
Concatenating those paths with a fixed path in $(G-v)-V(B-b)$ from $b$ to $y$ gives $\lfloor (k-1)/2\rfloor$ paths in $G$ from $x$ to $y$ with the length condition.
Therefore, $\lvert V(B-b)\cap \{x,y\} \rvert \in \{0,2\}$.

Since $G-v$ is not 2-connected, there exists an end-block $B'$ of $G-v$ with $V(B'-b')\cap \{x,y\}=\emptyset$, where $b'$ is the cut-vertex of $G-v$ contained in $B'$.
It follows that $N_G(B'-b')=\{b',v\}$.
Since $G$ is 2-connected, $G$ has two disjoint paths $P_1,P_2$ from $\{x,y\}$ to $\{b',v\}$ and internally disjoint from $B$.
Without loss of generality, we may assume that $P_1$ is from $x$ to $b'$ and $P_2$ is from $y$ to $v$.
Let $u$ be a vertex in $B'-b'$ adjacent to $v$ in $G$.
By Lemma \ref{path general graph}, $B'$ has $\lfloor (k-1)/2\rfloor$ paths $R_1,R_2,...,R_{\lfloor (k-1)/2\rfloor}$ from $b'$ to $u$ with the length condition.
Then $P_1\cup R_i\cup uv\cup P_2$ for all $i$ are $\lfloor (k-1)/2\rfloor$ paths in $G$ from $x$ to $y$ with the length condition.
This proves the lemma.
\end{pf}

\section{Cycles with the length condition} \label{sec:cycle}
\noindent In this section, we consider cycles with the length condition.
We first prove Theorem \ref{intro:cycle bipartite} in bipartite graphs.
We restate Theorem \ref{intro:cycle bipartite} here for the convenience of readers.

\medskip

\noindent {\bf Theorem \ref{intro:cycle bipartite}.}
{\em
Let $G$ be a bipartite graph and $v$ a vertex of $G$.
If every vertex of $G$ other than $v$ has degree at least $k+1$, then $G$ contains $k$ cycles with the length condition.
}

\medskip

\begin{pf}
Since there is nothing to prove when $k=0$, we may assume that $k \geq 1$.
We define a 2-connected end-block $H$ of $G$ and an edge $xy\in E(H)$ as following.
If $G$ is 2-connected, define $H=G$, $x=v$ and $y$ to be any neighbor of $x$ in $G$; if $G$ is not 2-connected, then define $H$ to be an end-block of $G$ such that $v\notin V(H-h)$, where $h$ is the cut-vertex of $G$ contained in $H$, and define $x= h$ and $y$ to be any neighbor of $x$ in $H$.
In either case, we see that every vertex of $H$ other than $x$ has degree at least $k+1$,
and thus $H$ is 2-connected bipartite with at least three vertices.
By Theorem \ref{thm:bipartite_main}, $H$ has $k$ paths from $x$ to $y$ with the length condition.
Note that each path has length at least two and thus does not contain the edge $xy$.
By adding the edge $xy$, we then obtain $k$ cycles in $H$ (and hence in $G$) with the length condition.
\end{pf}

\bigskip

\noindent {\bf Remark.} From the above proof,  it is easy to see that if $G$ is 2-connected bipartite with minimum degree at least $k+1$, then for every edge $e$ of $G$, there are $k$ cycles in $G$ with the length condition, and all of those cycles contain $e$.

\medskip

We then draw our attention to general graphs and prove Theorem \ref{intro:even/odd cycle in general graph}, which provides optimal bounds for cycles of consecutive even lengths as well as consecutive odd lengths.

\medskip

\noindent {\bf Theorem \ref{intro:even/odd cycle in general graph}.}
{\em
If the minimum degree of graph $G$ is at least $k+1$, then $G$ contains $\lfloor k/2\rfloor$ cycles with consecutive even lengths.
Furthermore, if $G$ is 2-connected and non-bipartite, then $G$ contains $\lfloor k/2\rfloor$ cycles with consecutive odd lengths.
}

\medskip

\begin{pf}
We may assume that $k\ge 2$, as the case $k=1$ is trivial.
Let $G'$ be a spanning bipartite subgraph of $G$ with the maximum number of edges,
and let $(A,B)$ be the bipartition of $G'$.
If $G'$ contains a vertex, say $v\in A$, of degree at most $\lfloor k/2 \rfloor$ in $G'$, then $(A-v,B\cup \{v\})$ will induce a bipartite subgraph of $G$ with more edges than $G'$, a contradiction.
So $G'$ has minimum degree at least $\lfloor k/2 \rfloor+1$.
By Theorem \ref{intro:cycle bipartite}, $G'$ (and hence $G$) contains $\lfloor k/2 \rfloor$ cycles with the length condition.
Note that each of those cycle has even length as $G'$ is bipartite.

Now we assume that $G$ is 2-connected and non-bipartite additionally.
Note that by the maximality, $G'$ is connected and bipartite with minimum degree at least $\lfloor k/2 \rfloor+1$.
Suppose that $G'$ is 2-connected.
Since $G$ is non-bipartite, there exist two vertices $x,y$ such that $xy\in E(G)-E(G')$.
So both $x,y$ are in the same part of the bipartition $(A,B)$.
By Theorem \ref{thm:bipartite_main}, $G'$ has $\lfloor k/2 \rfloor$ paths from $x$ to $y$ with the length condition.
Since both of $x,y$ are in the same part in the bipartition, each of these paths of $G'$ has even length. By concatenating these paths with the edge $xy$, we obtain $\lfloor k/2 \rfloor$ cycles in $G$ with consecutive odd lengths.
Hence, $G'$ is not 2-connected.
Let $H$ be an end-block of $G'$ and $h$ be the cut-vertex of $G'$ contained in $H$.
Every vertex of $H$ other than $h$ has degree at least $\lfloor k/2 \rfloor+1\ge 2$, so $H$ is 2-connected.
Since $G$ is 2-connected, there exist $z\in V(H-h)$ and $w\in V(G)-V(H)$ such that $zw\in E(G)-E(G')$.
By Theorem \ref{thm:bipartite_main}, $H$ has $\lfloor k/2 \rfloor$ paths from $z$ to $h$ with the length condition,
which, together with a fixed path in $G'-V(H-h)$ from $h$ to $w$, give $\lfloor k/2 \rfloor$ paths in $G'$ from $z$ to $w$ with the length condition.
As $zw\in E(G)-E(G')$, $z$ and $w$ are in the same part in the bipartition, so each of those mentioned paths in $G'$ from $z$ to $w$ has even length.
By concatenating these paths with the edge $zw$, we obtain $\lfloor k/2 \rfloor$ cycles in $G$ with consecutive odd lengths. This proves the theorem.
\end{pf}

\bigskip

\noindent {\bf Remark.} In fact, we can obtain $\lfloor k/2\rfloor$ cycles in $G$ with consecutive even lengths under a weaker condition that all vertices of $G$, but one, have degree at least $k+1$. On the other hand, we do not know if this weaker condition can guarantee the existence of $\lfloor k/2\rfloor$ consecutive odd cycles in Theorem \ref{intro:even/odd cycle in general graph}.

\medskip

As an immediate corollary of Theorem \ref{intro:even/odd cycle in general graph},
we can derive Theorem \ref{intro:cycles mod even k}, which proves Conjectures \ref{conj:Thom1} and \ref{conj:Thom2} when $k$ is even.

\medskip

\noindent {\bf Theorem \ref{intro:cycles mod even k}.}
{\em
Let $k$ be a positive even integer.
If the minimum degree of graph $G$ is at least $k+1$, then $G$ contains cycles of all even lengths modulo $k$.
Furthermore, if $G$ is 2-connected and non-bipartite, then $G$ contains cycles of all lengths modulo $k$.
}

\medskip

Theorem \ref{intro:even/odd cycle in general graph} also can be used to prove Theorem \ref{intro:ce(G)+co(G)}, which gives a tight relation between chromatic number and the number of cycles with the length condition.

\medskip

\noindent {\bf Theorem \ref{intro:ce(G)+co(G)}.}
{\em
For every graphs $G$, $\chi(G)\le 2 \min\{ce(G),co(G)\}+3$.
}

\medskip

\begin{pf}
We may assume that $\chi(G) \geq 3$, otherwise the theorem is easy.
Let $G'$ be a $\chi(G)$-critical subgraph of $G$.
Since $G'$ is $\chi(G)$-critical, $G'$ is 2-connected non-bipartite and $G'$ has minimum degree at least $\chi(G)-1$.
By Theorem \ref{intro:even/odd cycle in general graph}, $G'$ contains $\lfloor \chi(G)/2 \rfloor -1$ cycles with consecutive even lengths and contains $\lfloor \chi(G)/2 \rfloor -1$ cycles with consecutive odd lengths.
Hence $\min\{ce(G'),co(G')\} \geq \lfloor \chi(G)/2 \rfloor -1$.
As every cycle in $G'$ is a cycle in $G$, $\min\{ce(G),co(G)\} \geq \min\{ce(G'),co(G')\}\geq \lfloor \chi(G)/2 \rfloor -1 \geq (\chi(G)-1)/2-1$.
This proves the theorem.
\end{pf}

\bigskip

We conclude this section by proving a lemma about cycles with the length condition.

\begin{lemma} \label{2 but not 3-connected}
Let $G$ be a 2-connected but not 3-connected graph.
If the minimum degree of $G$ is at least $k+1$, then $G$ contains $2\lfloor k/2\rfloor-1$ cycles satisfying the length condition.
Furthermore, if $G$ is bipartite, then $G$ contains $2k-1$ cycles satisfying the length condition.
\end{lemma}

\begin{pf}
If $G$ is bipartite, let $t=k$; otherwise, let $t=\lfloor k/2\rfloor$.
Hence, by Theorem \ref{main bipartite} and Lemma \ref{path general graph}, for any subgraph $G'$ of $G$, if $(G',x,y)$ is 2-connected with minimum degree at least $t+1$, then $G'$ has $t$ paths from $x$ to $y$ with the length condition.
We shall prove that $G$ contains $2t-1$ cycles satisfying the length condition.

Since $G$ is 2-connected but not 3-connected,
there exists a separation $(A,B)$ of $G$ of order two.
Let $A\cap B=\{u,v\}$.
One can easily verify that each of $(G[A],u,v)$ and $(G[B],u,v)$ is a 2-connected rooted graph with minimum degree at least $k+1$.
Therefore, $G[A]$ has $t$ paths $P_1,P_2,...,P_t$ from $u$ to $v$ with the length condition, and $G[B]$ has $t$ paths $R_1,R_2,...,R_t$ from $u$ to $v$ with the length condition.
Then $P_i\cup R_j$ for all $1\leq i,j\leq t$ are $2t-1$ cycles satisfying the length condition.
\end{pf}

\section{Consecutive cycles} \label{sec:consecutive_cycles}

We say that a cycle $C$ in a connected graph $G$ is {\it non-separating} if $G-V(C)$ is connected.
The following lemma studies some property of non-separating odd cycle, which is a slight extension of \cite[Lemma 3.4]{Fan02}.

\begin{lemma} \label{lem:better non-separating cycle}
Let $G$ be a graph with minimum degree at least four.
If $G$ contains a non-separating induced odd cycle, then $G$ contains a non-separating induced odd cycle $C$, denoted by $v_0v_1...v_{2s}v_0$, such that either
	\begin{enumerate}
		\item $C$ is a triangle, or
		\item for every non-cut-vertex $v$ of $G-V(C)$, $\lvert N_G(v)\cap V(C) \rvert \le 2$, and the equality holds if and only if $N_G(v)\cap V(C)=\{v_i,v_{i+2}\}$ for some $i$, where the indices are taken under the additive group $\mathbb{Z}_{2s+1}$.
	\end{enumerate}
\end{lemma}

\begin{pf}
Let $C$ be a shortest non-separating induced odd cycle in $G$.
We denote $C=v_0v_1...v_{2s}v_0$.
Let $v$ be a non-cut-vertex of $G-V(C)$, and let $N_G(v)\cap V(C)=\{v_{i_1},...,v_{i_t}\}$ for some integers $i_1,...,i_t$ with $0\le i_1< ...< i_t\le 2s$.
Without loss of generality, we may assume that $i_1=0$.
For every $1\le j\le t$, let $C_j$ be the cycle $vv_{i_j}v_{i_j+1}...v_{i_{j+1}}v$.
Since the minimum degree of $G$ is at least four, every vertex in $C$ has at least one neighbor in $G-v-V(C)$, implying that $C_j$ is non-separating.
If $i_{j-1}=i_j+1$ for some $j$, then clearly $C_j$ is a non-separating triangle and hence $C$ is a triangle by the minimality.
So we may assume that $i_{j+1}-i_j \geq 2$, for each $j$ with $1 \leq j \leq t-1$, and $(2s+1)-i_t \geq 2$.
If $t\ge 3$, then for some $j$ the length of $C_j$ is odd and less than the length of $C$.
But $C_j$ is induced and non-separating, a contradiction to the minimality of $|V(C)|$.
So $t \leq 2$.
When $t=2$, by the minimality of $|V(C)|$, the unique even path in $C$ from $v_{i_1}$ to $v_{i_2}$ has to be of length two. This completes the proof.
\end{pf}

\begin{theorem} \label{nonseparating cycle->consecutive cycles}
Let $G$ be a 2-connected graph containing a non-separating induced odd cycle.
If the minimum degree of $G$ is at least $k+1$, then $G$ contains $2\lfloor \frac{k-1}{2}\rfloor$ cycles with consecutive lengths.
\end{theorem}

\begin{pf}
The theorem is obvious when $k\leq 2$.
So we may assume that $k\geq 3$.
By Lemma \ref{lem:better non-separating cycle}, there exists a non-separating induced odd cycle $C=v_0v_1...v_{2s}v_0$ in $G$ satisfying the conclusions of Lemma \ref{lem:better non-separating cycle}. 
Throughout this proof, the subscripts will be taken in the additive group $\mathbb{Z}_{2s+1}$.

\medskip

{\bf Claim 1:} $s\ge 2$ and hence $C$ is not a triangle. 

{\bf Proof of Claim 1.}
Suppose to the contrary that $C$ is a non-separating triangle $abca$.
Let $G' = (G-c)-\{ab\}$.
Since $G$ is 2-connected, $(G',a,b)$ is a 2-connected rooted graph with minimum degree at least $k$.
By Lemma \ref{path general graph}, $G'$ contains a sequence of $\lfloor (k-1)/2 \rfloor$ paths $P_1,...,P_{\lfloor (k-1)/2 \rfloor}$ from $a$ to $b$ satisfying the length condition.
Then $P_i\cup ba$ and $P_i\cup bca$ for all $1\le i\le \lfloor (k-1)/2 \rfloor$ are $2\lfloor (k-1)/2 \rfloor$ cycles in $G$ with consecutive lengths.
$\Box$

\medskip

So every non-cut-vertex $v$ of $G-V(C)$ has $d_{G-V(C)}(v)\ge k-1$.
Note that $s$ is a generator of the additive group ${\mathbb Z}_{2s+1}$.
For each $0\leq i \leq 2s$, let $v_i'= v_{i+s}$ and $v_i''=v_{i+s+1}$.
For any two vertices $v_i,v_j$ in $C$, denote $C'_{i,j}$ and $C''_{i,j}$ to be the shorter and longer paths in $C$ from $v_i$ to $v_j$, respectively.

\medskip

{\bf Claim 2:} $G-V(C)$ is not 2-connected.

\noindent{\bf Proof of Claim 2.}
Suppose to the contrary that $G-V(C)$ is 2-connected.
First assume that every vertex of $G-V(C)$ is adjacent in $G$ to at most one vertex of $C$. Then every vertex $v\in V(G-C)$ has $d_{G-V(C)}(v)\ge k$.
There exist distinct vertices $x,y\in V(G-C)$ such that $xv_0,yv_s\in E(G)$.
By Lemma \ref{path general graph}, $G-V(C)$ contains $\lfloor (k-1)/2 \rfloor$
paths $Q_1,...,Q_{\lfloor (k-1)/2 \rfloor}$ from $x$ to $y$ with the length condition.
Note that $C'_{0,s}$ and $C''_{0,s}$ are two paths from $v_0$ to $v_s$ of lengths $s,s+1$, respectively.
So $v_0x\cup Q_i\cup yv_s\cup C'_{0,s}$ and $v_0x\cup Q_i\cup yv_s\cup \cup C''_{0,s}$ for all $1\le i\le \lfloor (k-1)/2 \rfloor$ are $2\lfloor (k-1)/2 \rfloor$ cycles in $G$ with consecutive lengths.

Hence we may assume that there exists some $u\in V(G-C)$ adjacent to two vertices of $C$ in $G$. Without loss of generality, let $N_G(u)\cap V(C)=\{v_1,v_{2s}\}$,
and let $w\in V(G-C)$ such that $wv_s\in E(G)$.
Since $G-V(C)$ is 2-connected with minimum degree at least $k-1$,
by Lemma \ref{path general graph}, $G-V(C)$ contains a sequence of $\lfloor (k-2)/2 \rfloor$ paths $R_1,...,R_{\lfloor (k-2)/2 \rfloor}$ from $u$ to $w$ with the length condition.
Observe that $C'_{1,s}$ and $C'_{s,2s}$ are two paths of lengths $s-1$ and $s$, respectively and internally disjoint from $\{v_0,v_1,v_{2s}\}$.
Thus, $v_1u\cup R_i\cup wv_s\cup C'_{1,s}$ and $v_{2s}u\cup R_i\cup wv_s\cup C'_{s,2s}$ for all $1\le i\le \lfloor (k-2)/2 \rfloor$
together with $v_1v_0v_{2s}u\cup R_{\lfloor (k-2)/2 \rfloor}\cup wv_s\cup C'_{1,s}$ and $v_{2s}v_0v_1u\cup R_{\lfloor (k-2)/2 \rfloor}\cup wv_s\cup C'_{s,2s}$ give $2\lfloor k/2 \rfloor$ cycles in $G$ with consecutive lengths.
$\Box$

\medskip

Let $B$ be an end-block of $G-V(C)$ and $b$ the cut-vertex of $G-V(C)$ contained in $B$.
Every vertex in $B-b$ has degree at least $k-1\ge 2$ in $B$, and so $B$ is 2-connected.

\medskip

{\bf Claim 3:} There exists $x\in V(B-b)$ such that $N_G(x)\cap V(C)=\{v_{j-1},v_{j+1}\}$ for some $j$.

\noindent{\bf Proof of Claim 3.}
Suppose not that every vertex in $B-b$ is adjacent in $G$ to at most one vertex of $C$.
Then every vertex $v\in V(B-b)$ has $d_B(v)\ge k$.
If there exist $x\in V(B-b)$ and $y\in V(G-C)-V(B-b)$ such that $v_jx, v_j'y\in E(G)$ for some $j$,
then by Lemma \ref{path general graph}, $B$ contains $\lfloor (k-1)/2 \rfloor$ paths
$P_1,...,P_{\lfloor (k-1)/2 \rfloor}$ from $x$ to $b$ with the length condition.
Let $P$ be a path in $G-V(C)-V(B-b)$ from $b$ to $y$.
Also note that $C'_{j,j+s}$ and $C''_{j,j+s}$ are two paths in $C$ from $v_j$ to $v_j'$ of lengths $s$ and $s+1$, respectively.
Then, $v_jx\cup P_i\cup P\cup yv_j'\cup C'_{j,j+s}$ and $v_jx\cup P_i\cup P\cup yv_j'\cup C''_{j,j+s}$ for all $1\le i\le \lfloor (k-1)/2 \rfloor$ are $2\lfloor (k-1)/2 \rfloor$ cycles in $G$ with consecutive lengths.
Hence, we may assume that if $v_j$ is adjacent to $V(B-b)$, then $N_G(v_j') \cap V(G-C) \subseteq V(B-b)$.
There is some vertex of $C$ adjacent in $G$ to $V(B-b)$, and $s$ is a generator of ${\mathbb Z}_{2s+1}$, so we derive that $N_G(C)\subseteq V(B-b)$.
This implies that $b$ is a cut-vertex of $G$, but $G$ is 2-connected, a contradiction.
$\Box$

\medskip

{\bf Claim 4:} $N_G(\{v'_j,v''_j\})\cap V(G-C)\subseteq V(B-b)$.

\noindent{\bf Proof of Claim 4.}
Suppose not, by symmetry we may assume that $v'_jy\in E(G)$ for some $y\in V(G-C)-V(B-b)$. Since every vertex in $B-b$ has degree at least $k-1$,
by Lemma \ref{path general graph}, $B$ contains $\lfloor (k-2)/2 \rfloor$ paths $Q_1,...,Q_{\lfloor (k-2)/2 \rfloor}$ from $x$ to $b$ with the length condition.
Let $Q$ be a fixed path in $G-V(C)-V(B-b)$ from $b$ to $y$.
Note that $C'_{j+1,j+s}, C'_{j-1,j+s}$ are two paths in $C$ from $v_j'$ to $v_{j+1}, v_{j-1}$ with lengths $s-1, s$, respectively and internally disjoint from $\{v_{j-1},v_j,v_{j+1}\}$.
Then, $v_{j+1}x\cup Q_i\cup Q\cup yv'_j\cup C'_{j+1,j+s}$ and
$v_{j-1}x\cup Q_i\cup Q\cup yv'_j\cup C'_{j-1,j+s}$ for all $1\le i\le \lfloor (k-2)/2 \rfloor$, together with $v_{j+1}v_jv_{j-1}x\cup Q_{\lfloor (k-2)/2 \rfloor}\cup Q\cup yv'_j\cup C'_{j+1,j+s}$ and $v_{j-1}v_jv_{j+1}x\cup Q_{\lfloor (k-2)/2 \rfloor}\cup Q\cup yv'_j\cup C'_{j-1,j+s}$, are $2\lfloor k/2 \rfloor$ cycles in $G$ with consecutive lengths.
$\Box$

\medskip

Since $d_{G-V(C)}(v'_j)\ge k-1\ge 2$, there exists $z\in V(B)-\{b,x\}$ adjacent to $v'_j$.
Every vertex of $B$ other than $b$ has degree at least $k-1$ in $B$.
By Lemma \ref{lem:path with one exceptional vertex}, $B$ has $\lfloor (k-3)/2 \rfloor$ paths $R_1,...,R_{\lfloor (k-3)/2 \rfloor}$ from $x$ to $z$ with the length condition.
Then, $v_{j+1}x\cup R_i\cup zv'_j\cup C'_{j+1,j+s}$ and $v_{j-1}x\cup R_i\cup zv'_j\cup C'_{j-1,j+s}$ for all $1\le i\le \lfloor (k-3)/2 \rfloor$,
together with $v_{j+1}v_jv_{j-1}x\cup R_{\lfloor (k-3)/2 \rfloor}\cup zv'_j\cup C'_{j+1,j+s}$ and $v_{j-1}v_jv_{j+1}x\cup R_{\lfloor (k-3)/2 \rfloor}\cup zv'_j\cup C'_{j-1,j+s}$, are $\lfloor (k-1)/2 \rfloor$ cycles in $G$ with consecutive lengths. This completes the proof of Theorem \ref{nonseparating cycle->consecutive cycles}.
\end{pf}

\bigskip

Now we are ready to prove Theorems \ref{intro:3-connected non-bipartite} and \ref{intro:2-con non-bipa cycles}.

\medskip

\noindent{\bf Theorem \ref{intro:3-connected non-bipartite}.}
{\em
If $G$ is a 3-connected non-bipartite graph with minimum degree at least $k+1$, then $G$ contains $2\lfloor \frac{k-1}{2}\rfloor$ cycles with consecutive lengths.
}

\medskip

\begin{pf}
It was proved by several groups (see \cite{TT81,BV98}) that every 3-connected non-bipartite graph contains a non-separating induced odd cycle.
This, together with Theorem \ref{nonseparating cycle->consecutive cycles}, immediately imply this theorem.
\end{pf}


\medskip

\noindent{\bf Theorem \ref{intro:2-con non-bipa cycles}.}
{\em
If $G$ is a 2-connected non-bipartite graph with minimum degree at least $k+3$,
then $G$ contains $k$ cycles with consecutive lengths or the length condition.
}

\medskip

\begin{pf}
If $G$ is 3-connected, then by Theorem \ref{intro:3-connected non-bipartite},
$G$ contains $2\lfloor (k+1)/2\rfloor\ge k$ cycles with consecutive lengths.
Otherwise $G$ is 2-connected but not 3-connected, by Lemma \ref{2 but not 3-connected}, $G$ contains $2\lfloor (k+2)/2\rfloor-1\ge k$ cycles with the length condition.
\end{pf}

\bigskip

From this result, we can prove Theorem \ref{intro:2-con non-bipa cycles modulo k} promptly.

\medskip

\noindent{\bf Theorem \ref{intro:2-con non-bipa cycles modulo k}.}
{\em
Let $k$ be a positive odd integer.
If $G$ is a 2-connected non-bipartite graph with minimum degree at least $k+3$,
then $G$ contains cycles of all lengths modulo $k$.
}

\medskip

\begin{pf}
By Theorem \ref{intro:2-con non-bipa cycles}, $G$ contains $k$ cycles with consecutive lengths or the length condition. Since $k$ is odd, in either case, the set of these cycle lengths intersect each of the residue classes modulo $k$.
\end{pf}

\bigskip

The following theorem will be used for proving Theorem \ref{intro:more general graph k cycles}.

\begin{theorem} \label{general graph k cycles}
Let $G$ be a 2-connected graph and $v$ a vertex of $G$.
If every vertex of $G$ other than $v$ has degree at least $k+4$, then $G$ contains $k$ cycles with consecutive lengths or with the length condition.
\end{theorem}

\begin{pf}
Let $G'=G-\{v\}$.
So $G'$ has minimum degree at least $k+3$.
If $G'$ is bipartite, then $G'$ contains $k+2$ cycles with the length condition by Theorem \ref{intro:cycle bipartite}.
So we may assume that $G'$ is non-bipartite.

If $G'$ is 2-connected, then by Theorem \ref{intro:2-con non-bipa cycles}, $G'$ contains
$k$ cycles with consecutive lengths or the length condition.
So we may assume that $G'$ is not 2-connected.
Note that the minimum degree of $G'$ is at least $k+3$, so every end-block of $G'$ is 2-connected.

Since $G$ is 2-connected, $G'$ contains two end-blocks $B_1,B_2$ such that for each $i \in \{1,2\}$, $B_i-b_i$ contains a vertex $v_i$ adjacent in $G$ to $v$, where $b_i$ is the cut-vertex of $G'$ contained in $B_i$.
By Lemma \ref{path general graph}, for each $i \in \{1,2\}$, $B_i$ contains $\lfloor (k+2)/2\rfloor$ paths $P_{i,1},...,P_{i,\lfloor (k+2)/2\rfloor}$ from $b_i$ to $v_i$ with the length condition.
Let $R$ be a path in $G'$ from $b_1$ to $b_2$ internally disjoint from $V(B_1)\cup V(B_2)$.
Then for $1\le j,j'\le \lfloor (k+2)/2\rfloor$, $P_{1,j}\cup R\cup P_{2,j'}\cup v_2vv_1$ are $2\lfloor (k+2)/2\rfloor-1\ge k$ cycles in $G$ with the length condition.
\end{pf}

\medskip

\noindent{\bf Theorem \ref{intro:more general graph k cycles}.}
{\em
If $G$ is a graph with minimum degree at least $k+4$,
then $G$ contains $k$ cycles with consecutive lengths or the length condition.
}

\medskip

\begin{pf}
Let $B$ be an end-block of $G$ and let $b$ be the cut-vertex of $G$ contained in $B$.
Every vertex of $B$ other than $b$ has minimum degree at least $k+4$ and hence $B$ is 2-connected.
By Theorem \ref{general graph k cycles}, $B$ (and hence $G$) contains $k$ cycles with consecutive lengths or with the length condition.
\end{pf}

\bigskip

It is straightforward to obtain Theorem \ref{intro:general cycles mod k} from Theorem \ref{intro:more general graph k cycles}.

\medskip

\noindent{\bf Theorem \ref{intro:general cycles mod k}.}
{\em
Let $k$ be a positive odd integer. If $G$ is a graph with minimum degree at least $k+4$,
then $G$ contains cycles of all lengths modulo $k$.
}

\medskip

\begin{pf}
By Theorem \ref{intro:more general graph k cycles}, $G$ contains $k$ cycles with consecutive lengths or the length condition. Since $k$ is odd, in either case, the set of these cycle lengths intersect each of the residue classes modulo $k$.
\end{pf}

\bigskip

Lastly, we derive Theorem \ref{intro:c(G)} from Theorem \ref{nonseparating cycle->consecutive cycles}.

\medskip

\noindent{\bf Theorem \ref{intro:c(G)}.}
{\em
For every graphs $G$, $\chi(G)\le c(G)+4$.
}

\medskip

\begin{pf}
Suppose to the contrary that there exists a graph $G$ with $\chi(G)\ge c(G)+5$.
Let $G'$ be a $\chi(G)$-critical subgraph of $G$.
Note that $G'$ is 2-connected and has minimum degree at least $\chi(G)-1 \geq c(G)+4$.
A result of Krusenstjerna-Hafstr{\o}m and Toft (see \cite{KHT}, Theorem 4) states that every 4-critical graph contains a non-separating induced odd cycle, but in fact their proof also works for $k$-critical graph for every $k \geq 4$.
(We direct interested readers to the original proof in \cite{KHT}.)
Thus, $G'$ also contains a non-separating induced odd cycle.
By Theorem \ref{nonseparating cycle->consecutive cycles}, $G'$ contains $2\lfloor \frac{c(G)+2}{2}\rfloor\ge c(G)+1$ consecutive cycles.
However, every cycle in $G'$ is a cycle in $G'$, so $c(G) \geq c(G')\geq c(G)+1$, a contradiction.
This completes the proof.
\end{pf}

\section{Concluding remarks}\label{sec:remarks}
\noindent In this paper, we have obtained several tight or nearly tight results on the relation between cycle lengths and minimum degree.
It will be interesting if one can close the gap between our results and the best possible upper bounds, such as in Theorems \ref{intro:3-connected non-bipartite} and \ref{intro:2-con non-bipa cycles}.
A good starting point may be the following strengthening of Theorem \ref{intro:even/odd cycle in general graph}.

\begin{conjecture}
If $G$ is a 2-connected non-bipartite graph with minimum degree at least $k+1$, then $G$ contains $\lceil k/2\rceil$ cycles with consecutive odd lengths.
\end{conjecture}

\noindent If it is true, then one can prove $\chi(G)\le 2 co(G)+2$ as in Theorem \ref{intro:ce(G)+co(G)}.

In Theorem \ref{intro:more general graph k cycles}, we prove that every graph $G$ with $\delta(G)\ge k+4$ contains $k$ cycles with consecutive lengths or the length condition.
The following examples show that the bound $\delta(G)\ge k+4$ is tight up to the constant term: the complete graph $K_{k+2}$ has precisely $k$ cycles of consecutive lengths $3,4,...,k+2$, while for every $n\ge k+1$ the complete bipartite graph $K_{k+1,n}$ has precisely $k$ cycles of consecutive even lengths $4,6,..,2k+2$.
All such graphs have minimum degree $k+1$, and thus we conjecture that $\delta(G)\ge k+1$ is optimal.

\begin{conjecture}
Every graph with minimum degree at least $k+1$ contains $k$ cycles with consecutive lengths or the length condition.
\end{conjecture}

\noindent If true, this would imply both Conjectures \ref{conj:Thom1} and \ref{conj:Thom2} when $k$ is odd, and thus, together with Theorems \ref{intro:cycles mod even k}, imply these conjectures in full generality.

Our results show that if a graph $G$ has $\delta(G) \geq k+4$ (and satisfies some necessary conditions), then $G$ contains cycles of all lengths modulo $k$.
This is tight up to the constant term.
However, for fixed integer $m$, we know very little about the least function $f(m,k)$ such that every graph $G$ with $\delta(G) \geq f(m,k)$ contains a cycle of length $m$ modulo $k$. (If $k$ is even and $m$ is odd, then one has to restrict to 2-connected non-bipartite graphs $G$ here.)
A conjecture of Dean (see \cite{DLS93}) considered the case when $m=0$, which asserted that every $k$-connected graph contains a cycle of length $0$ modulo $k$.
Note that this (if true) is best possible for odd $k$,
as for every $n\ge k-1$, $K_{k-1,n}$ is $(k-1)$-connected but has no cycles of length $0$ modulo $k$.
Dean's conjecture was confirmed for $k=3$ in \cite{CS} and $k=4$ in \cite{DLS93}.
Another interesting special case is $m=3$ (for the sake of convenience, let $k$ be odd).
So $f(3,k)$ becomes the least function such that every triangle-free graph $G$ with minimum degree $f(3,k)$ contains a cycle of length 3 modulo $k$.
We speculate that $f(3,k)=o(k)$. This may be related to the recent result of \cite{KSV}.

Despite much research has been done, the distribution of cycle lengths in graphs with large minimum degree is still mysterious and unclear.
We conclude this paper by mentioning a conjecture of Erd\H{o} and Gy\'arf\'as \cite{Erd93}:
every graph with minimum degree at least three contains a cycle of length a power of two.

\end{document}